\documentclass[11pt]{amsart}
\usepackage{amsmath,amsfonts}
\usepackage{graphicx}
\usepackage{color}

\newtheorem{thm}{Theorem}[section]
\newtheorem{cor}[thm]{Corollary}
\newtheorem{lem}[thm]{Lemma}
\newtheorem{fact}[thm]{Fact}

\newtheorem{rem}[thm]{Remark}
\newtheorem{ex}[thm]{Example}

\numberwithin{equation}{section}

\newcommand{\dist}{\operatorname{dist}}
\newcommand{\co}{\operatorname{cone}}

\begin{document}

\title[Approximation by smooth mappings with no critical points]
{Approximation by smooth functions with no critical points on
separable Banach spaces}
\author{D. Azagra and M. Jim\'{e}nez-Sevilla}
\date{October, 2005}

\thanks{\noindent D. Azagra was supported by a
Marie Curie Intra-European Fellowship of the European Community,
Human Resources and Mobility Programme under contract number MEIF
CT2003-500927.\\ M. Jim\'{e}nez-Sevilla was supported by a Fellowship
of the {Secretar\'{\i}}a de Estado de Universidades e Investigaci\'{o}n
(Ministerio de Educaci\'{o}n y Ciencia). }

\begin{abstract}
We characterize the class of separable Banach spaces $X$ such that
for every continuous function $f:X\to\mathbb{R}$ and for every
continuous function $\varepsilon:X\to\mathbb(0,+\infty)$ there
exists a $C^1$ smooth function $g:X\to\mathbb{R}$ for which
$|f(x)-g(x)|\leq\varepsilon(x)$ and $g'(x)\neq 0$ for all $x\in X$
(that is, $g$ has no critical points), as those Banach spaces $X$
with separable dual $X^*$. We also state sufficient conditions on
a separable Banach space so that the function $g$ can be taken to
be of class $C^p$, for $p=1,2,..., +\infty$. In particular, we
obtain the optimal order of smoothness of the approximating
functions with no critical points on the classical spaces
$\ell_p(\mathbb{N})$ and $L_p(\mathbb{R}^n)$. Some important
consequences of the above results are (1) the existence of {\em a
non-linear Hahn-Banach theorem} and (2) the smooth approximation
of closed sets, on the classes of spaces considered above.

\end{abstract}

\maketitle

\section[Introduction and main results]{Introduction and main results}

The Morse-Sard theorem \cite{Sard1, Sard2} states that if
$f:\mathbb{R}^{n}\longrightarrow \mathbb{R}^{m}$ is a $C^r$ smooth
function, with $r>\max\{n-m, 0\}$, and $C_{f}$ is the set of
critical points of $f$, then the set of critical values $f(C_{f})$
is of Lebesgue measure zero in $\mathbb{R}^{m}$. This result has
proven to be very valuable in a large number of areas, especially
in differential topology and analysis (see for instance
\cite{Hirsch, YomdinComte} and the references therein). Additional
geometric and analytical properties  of the set of critical values
in different versions of the Morse-Sard theorem, together with a
study on the sharpness of the hypothesis of the Morse-Sard
theorem, have been obtained in \cite{Bates1, Bates2, Bates3,
Bates4, Bates-Moreira, Moreira}.

For many important applications of the Morse-Sard theorem, it is
enough to know that any given continuous function can be uniformly
approximated by a smooth map whose set of critical values has
empty interior \cite{Hirsch, YomdinComte}. We refer to this as an
{\em approximate Morse-Sard theorem}. The same type of
approximation could prove key to the study of related problems in
the infinite-dimensional domain.

In this paper, we will prove the strongest version of an
approximate Morse-Sard theorem that one can expect to be true for
a general infinite-dimensional separable Banach space, namely that
every continuous function $f:X\longrightarrow \mathbb R$, where
$X$ is an infinite-dimensional Banach space $X$ with separable
dual $X^*$, can be uniformly approximated by a $C^1$ smooth
function $g:X\longrightarrow \mathbb R$ which does not have any
critical point. In some cases where more information about the
structure of the Banach space $X$ is known, we will  extend our
result to higher order of differentiability, $C^p$ ($p>1$).

Our result will also allow us to demonstrate two important
corollaries. The first one is the existence of {\em a non-linear
Hahn-Banach theorem} which shows that two disjoint closed subsets
in $X$ can be separated by a 1-codimensional $C^p$ smooth manifold
of $X$ (which is the set of zeros of a $C^p$ smooth function with
no critical points on $X$). The second one states that every
closed subset of $X$ can be approximated by $C^p$ smooth open
subsets of $X$.

To put our work in context, let us briefly review some of the work
established for the infinite-dimensional version of the Morse-Sard
theorem. Smale \cite{Smale} proved that if $X$ and $Y$ are
separable connected smooth manifolds modelled on Banach spaces and
$f:X\longrightarrow Y$ is a $C^r$ Fredholm map then $f(C_{f})$ is
of first Baire category and, in particular, $f(C_{f})$ has no
interior points provided that $r>\max\{\textrm{index}(df(x)), 0\}$
for all $x\in X$. Here, index($df(x)$) stands for the index of the
Fredholm operator $df(x)$, that is, the difference between the
dimension of the kernel of $df(x)$ and the codimension of the
image of $df(x)$, which are both finite. These assumptions are
very strong as they impose that when $X$ is infinite-dimensional
then $Y$ is necessarily infinite-dimensional too (in other words,
there is no Fredholm map $f:X\longrightarrow\mathbb{R}$). In fact,
as Kupka proved in \cite{Kupka}, there are $C^\infty$ smooth
functions $f:\ell_2\longrightarrow\mathbb{R}$ (where $\ell_2$ is
the separable Hilbert space) such that their sets of critical
values $f(C_{f})$ contain intervals and hence have non-empty
interiors and positive Lebesgue measure. Bates and Moreira
\cite{Bates-Moreira, Moreira} showed that this function $f$ can
even be taken to be a polynomial of degree three. Azagra and
Cepedello-Boiso \cite{AC} have shown that every continuous mapping
from the separable Hilbert space into $\mathbb{R}^{m}$ can be
uniformly approximated by $C^\infty$ smooth mappings with no
critical points. Unfortunately, since the core of their proof
requires the use of the special properties of the Hilbertian norm,
this  cannot be extended to non-Hilbertian Banach spaces. P.
H\'{a}jek and M. Johanis \cite{HJ} established the same kind of
result in the case when $X$ is a separable Banach space which
contains $c_0$ and admits a $C^p$-smooth bump function. In this
case, the approximating functions are of class $C^p$, $p=1, 2,...,
\infty$. This method is based on the result that the range of the
derivative of a $C^2$ smooth function from $c_0$ to $\mathbb R$ is
a countable union of compact sets \cite{Hajek}. However, as the
authors noted, their method is not applicable when the space $X$
has the Radon-Nikod\'{y}m property (e.g., when $X$ is reflexive),
which leaves out all the classical Banach spaces $\ell_p$ and
$L_{p}(\mathbb{R}^{n})$ for $1<p<\infty$.

As stated above, we prove that for any infinite-dimensional Banach
space $X$ with a separable dual $X^*$, the set of $C^1$ smooth,
real-valued functions with no critical points is uniformly dense
in the space of all continuous, real-valued functions on $X$. This
solves completely the problem of the approximation on separable
Banach spaces by smooth, real-valued functions with no critical
points when the order of smoothness of the approximating functions
is one. Hence, we obtain the following characterization. For a
separable Banach space $X$, the following are equivalent: (i)
$X^*$ is separable, and (ii) the set of $C^1$ smooth, real-valued
functions on $X$ with no critical points is uniformly dense in the
space of all continuous, real-valued functions on $X$.

This result can be included in our main theorem which also applies
to higher order of differentiability. Before stating our main
theorem, recall that a norm $||\cdot||$ in a Banach space $X$ is
LUR (locally uniformly rotund \cite{DGZ}) if $\lim_n||x_n-x||=0$
whenever the sequence $\{x_n\}_n$ and the point $x$ are included
in the unit sphere of the norm $||\cdot||$ and
$\lim_n||x_n+x||=2$. A norm $||\cdot||$ in $X$ is $C^p$ smooth if
it is $C^p$ smooth in $X\setminus\{0\}$.

\begin{thm}\label{approximation theorem}
Let $X$ be an infinite dimensional separable Banach space $X$ with
a LUR and $C^p$ smooth norm $||\cdot||$, where $p\in \mathbb N\cup
\{\infty\}$. Then, for every continuous  mapping
$f:X\longrightarrow\mathbb{R}$ and for every continuous function
$\varepsilon:X\longrightarrow (0,\infty)$, there exists a $C^p$
smooth mapping $g:X\longrightarrow\mathbb{R}$ such that
$|f(x)-g(x)|\leq\varepsilon(x)$ for all $x\in X$ and $g$ has no
critical points.
\end{thm}

Our proof involves: {\em i)} a special construction of carefully
perturbed partitions of unity in an open subset of the unit sphere
of the Banach space \, $Y=X\oplus\mathbb R$ \, by means of a
sequence of linear functionals in $Y^*$, {\em ii)} the study and
use of the properties of the range of the derivative of the norm
in $Y$, $Y^*$ and their finite dimensional subspaces (Lemmas
\ref{case N} and \ref{vectoradicional} below), and {\em iii)} the
use of $C^p$ deleting diffeomorphisms from $X$ onto $X\setminus
O$, where $O$ is a bounded, closed, convex subset of $X$.

The following example gives the optimal order of smoothness of the
approximation functions with no critical points for
$\ell_p(\mathbb N)$ and $L_{p}(\mathbb{R}^{n})$.

\begin{ex}\label{example}
It follows immediately from Theorem \ref{approximation theorem}
that one can approximate every continuous, real-valued  function
on $\ell_{p}(\mathbb{N})$ and $L_{p}(\mathbb{R}^{n})$ \
($1<p<\infty$) \ with $C^{\overline{p}}$ smooth, real-valued
functions with no critical points, where \ $\overline{p}=[p]$ \ if
$p$ is not an integer, \ $\overline{p}=p-1$ \ if $p$ is an odd
integer, and \ $\overline{p}=\infty$ \  if $p$ is an even integer.
Indeed, the standard norms of the classical separable Banach
spaces $\ell_{p}(\mathbb{N})$ and $L_{p}(\mathbb{R}^{n})$ are LUR
and $C^{\overline{p}}$ smooth \cite{DGZ}.
\end{ex}

\bigskip

Since every Banach space with separable dual admits an equivalent
LUR and $C^1$ smooth norm \cite{DGZ}, we immediately deduce from
Theorem \ref{approximation theorem} the announced characterization
of the property of approximation by $C^1$ smooth functions with no
critical points.

\begin{cor}\label{aproximacionC^1sinpuntoscriticos}
Let $X$ be a separable Banach space. The following are equivalent:

\begin{enumerate}
\item The dual space $X^*$ is separable, \item  for every
continuous mapping $f:X\longrightarrow\mathbb{R}$ and for every
continuous function $\varepsilon:X\longrightarrow (0,\infty)$,
there exists a $C^1$ smooth mapping $g:X\longrightarrow\mathbb{R}$
such that $|f(x)-g(x)|\leq\varepsilon(x)$ and $g$ has no critical
points.
\end{enumerate}
\end{cor}

\medskip

Next, we establish  a similar statement for higher order
smoothness on separable Banach spaces with a $C^p$ smooth bump
function ($p\ge 2$) and unconditional basis. We combine Theorem
\ref{approximation theorem} and the results on fine approximation
given in \cite{AFGJL}, to obtain the optimal order of smoothness
of the approximating functions with no critical points on a large
class within the Banach spaces with separable dual. In particular
the following Corollary applies even when the space $X$ lacks a
norm which is simultaneously LUR and $C^2$ smooth.

\medskip

\begin{cor} \label{fineaprox} Let $X$ be a  separable Banach space with unconditional basis.
Assume that $X$ has a $C^p$ smooth Lipschitz bump function, \
where $p\in \mathbb N\cup \{\infty\}$. Then, for every continuous
mapping $f:X\longrightarrow\mathbb{R}$ and for every continuous
function $\varepsilon:X\longrightarrow (0,\infty)$, there exists a
$C^p$ smooth mapping \ $g:X\longrightarrow\mathbb{R}$ such that
$|f(x)-g(x)|\leq\varepsilon(x)$ and $g$ has no critical points.
\end{cor}
\begin{proof}
Since $X$ is separable and admits a $C^p$ smooth bump function,
the dual space $X^*$ is separable. Thus we obtain, from Corollary
\ref{aproximacionC^1sinpuntoscriticos}, a $C^1$ smooth function
$h:X\longrightarrow \mathbb R$ such that \ $h'(x)\neq 0$ \ and \
$|f(x)-h(x)|<\frac{\varepsilon(x)}{2}$ \ for every  \ $x\in X$.
Let us denote by $||\cdot||$ the dual norm on $X^*$. Define the
continuous function \ $\overline{\varepsilon}:X\longrightarrow
(0,\infty)$, \
$\overline{\varepsilon}(x)=\frac12\min\{\varepsilon(x),\
||h'(x)||\}$, for $x\in X$. Now, by the main result of
\cite{AFGJL}, there is a $C^p$ smooth function
$g:X\longrightarrow\mathbb{R}$ such that \
$|h(x)-g(x)|<\overline{\varepsilon}(x)$ \ and \
$||h'(x)-g'(x)||<\overline{\varepsilon}(x)$, \ for every $x\in X$.
The latter implies that \ $||h'(x)||-||g'(x)||< \frac12
||h'(x)||$, \ and therefore \ $0<\frac12||h'(x)||<||g'(x)||$ \ for
every $x\in X$. Hence, \ $g$ is a $C^p$ \ smooth function with no
critical points and \
$|f(x)-g(x)|<|f(x)-h(x)|+|h(x)-g(x)|<\varepsilon(x)$ \ for every
$x\in X$.
\end{proof}

The proof of the above corollary yields to the following remark.
\begin{rem}
Assume that a separable Banach space $X$ satisfies  the $C^1$-fine
approximation property by $C^p$ smooth, real-valued functions,
i.e., for every $C^1$ smooth function $f:X\longrightarrow \mathbb
R$ and every continuous function $\varepsilon:X\longrightarrow
\mathbb (0,\infty)$ there is a $C^p$ smooth function
$h:X\longrightarrow \mathbb R$ such that $|f(x)-h(x)|\le
\varepsilon(x)$ and $|f'(x)-h'(x)|\le \varepsilon(x)$, for every
$x\in X$. Then, the conclusion of Corollary \ref{fineaprox} holds.
\end{rem}

\medskip

Furthermore, our results allow us to make the following
conclusions.

\begin{rem}
{\em
\hfill{ }
\begin{enumerate}
\item All of the results presented above hold in the case when one
replaces $X$ with an open subset $U$ of $X$. Actually, the same
proof given in the section to follow (with obvious modifications)
can be used.

\item Whenever $X$ has the property that every continuous,
real-valued function on $X$ can be approximated by $C^p$ smooth,
real-valued functions with no critical points, one can deduce the
following Corollaries.
\end{enumerate}
}
\end{rem}

\begin{cor}[A nonlinear Hahn-Banach theorem]\label{nonlinear Hahn-Banach theorem}
Let $X$ be any of the Banach spaces considered in the above
results. Then, for every two disjoint closed subsets $C_1$, $C_2$
of $X$, there exists a $C^p$ smooth function
$\varphi:X\longrightarrow\mathbb{R}$ with no critical points, such
that the level set $M=\varphi^{-1}(0)$ is a $1$-codimensional
$C^p$ smooth submanifold of $X$ that separates $C_{1}$ and
$C_{2}$, in the following sense: define $U_{1}=\{x\in X :
\varphi(x)<0\}$ and $U_{2}=\{x\in X : \varphi(x)>0\}$, then $U_1$
and $U_2$ are disjoint $C^p$ smooth open sets of $M$ with common
boundary $\partial U_{1}=\partial U_{2}=M$, and such that
$C_{i}\subset U_{i}$ for $i=1, 2$.
\end{cor}

Recall that an open subset $U$ of $X$ is said to be $C^p$ smooth
provided its boundary $\partial U$ is a $C^p$ smooth
one-codimensional submanifold of $X$.

\begin{cor}[Smooth approximation of closed sets]\label{Smooth
approximation of closed sets} Every closed subset of any of the
Banach spaces $X$ considered above can be approximated by $C^p$
smooth open subsets of $X$ in the following sense: for every
closed set $C\subset X$ and every open set $W$ containing $C$
there is a $C^p$ smooth open set $U$ so that $C\subset U\subseteq
W$.
\end{cor}

\begin{cor}[Failure of Rolle's Theorem]\label{failure of Rolle's theorem}
For every open subset $U$ of any of the above Banach spaces $X$
there is a continuous function $f$ on $X$ whose support is the
closure of $U$, and such that $f$ is $C^p$ smooth on $U$ and yet
$f$ has no critical point in $U$.
\end{cor}

\bigskip

\section[Proof]{Proof of Theorem \ref{approximation theorem}}

Recall that a norm $N(\cdot)$, in a Banach space $E$, is (1) {\em
strictly convex} if the unit sphere of the norm $N(\cdot)$ does
not include any segment line. Equivalently, $N(\frac{x+y}{2})<1$
for every $x,y$ in the unit sphere with $x\not=y$; (2) {\em WUR
(weakly uniformly rotund)} if $\lim_n(x_n-y_n)=0$ in the weak
topology whenever the sequences $\{x_n\}_n$ and $\{y_n\}_n$ are
included in the unit sphere and $\lim_nN(x_n+y_n)=2$.

\medskip

We denote by $\mathbb R^*$ the set of non zero real numbers. We
also let $[u_1,...,u_n]$ stand for the linear span of the vectors
$u_1$,...,$u_n$. Let us denote by $S_{||\cdot||}$ and
$S_{||\cdot||^*}$ the unit sphere of a Banach space
$(Z,||\cdot||)$ and its dual $(Z^*,||\cdot||^*)$, respectively.

\medskip

The following two geometrical lemmas will be essential to the
proof of Theorem \ref{approximation theorem}.

\bigskip

\begin{lem}\label{case N}
Let $Z=[u_1,...,u_n]$ be a n-dimensional space ($n>1$) with a
differentiable norm $||\cdot||$ (G\^{a}teaux or Fr\'{e}chet
differentiable, since both notions coincide for convex functions
defined on finite dimensional spaces). Let us consider real
numbers \ $0<\alpha_i<1$,\ for $i=1,...,n-1$ and define
$\mathcal{R}$ as the subset of
 numbers $\alpha\in\mathbb R^*$ satisfying that

\begin{equation}\label{vacio}
\big\{T\in S_{||\cdot||^*}:\ T(u_i)=\alpha_i\,, \ i=1,...,n-1, \
T(u_n)=\alpha \bigr\}=\emptyset.
\end{equation}
 Then, the cardinal of $\mathbb R^*\setminus \mathcal{R}$ is at most two.
\end{lem}
\begin{proof}
First, assume that the set
\begin{equation*}
F=\{T\in S_{||\cdot||^*}: \
T(u_1)=\alpha_1,...,T(u_{n-1})=\alpha_{n-1} \}\end{equation*} is
non-empty (otherwise we have finished).

\begin{figure}
\centering
\begin{tabular}{ccc}
\includegraphics[width=4.5cm]{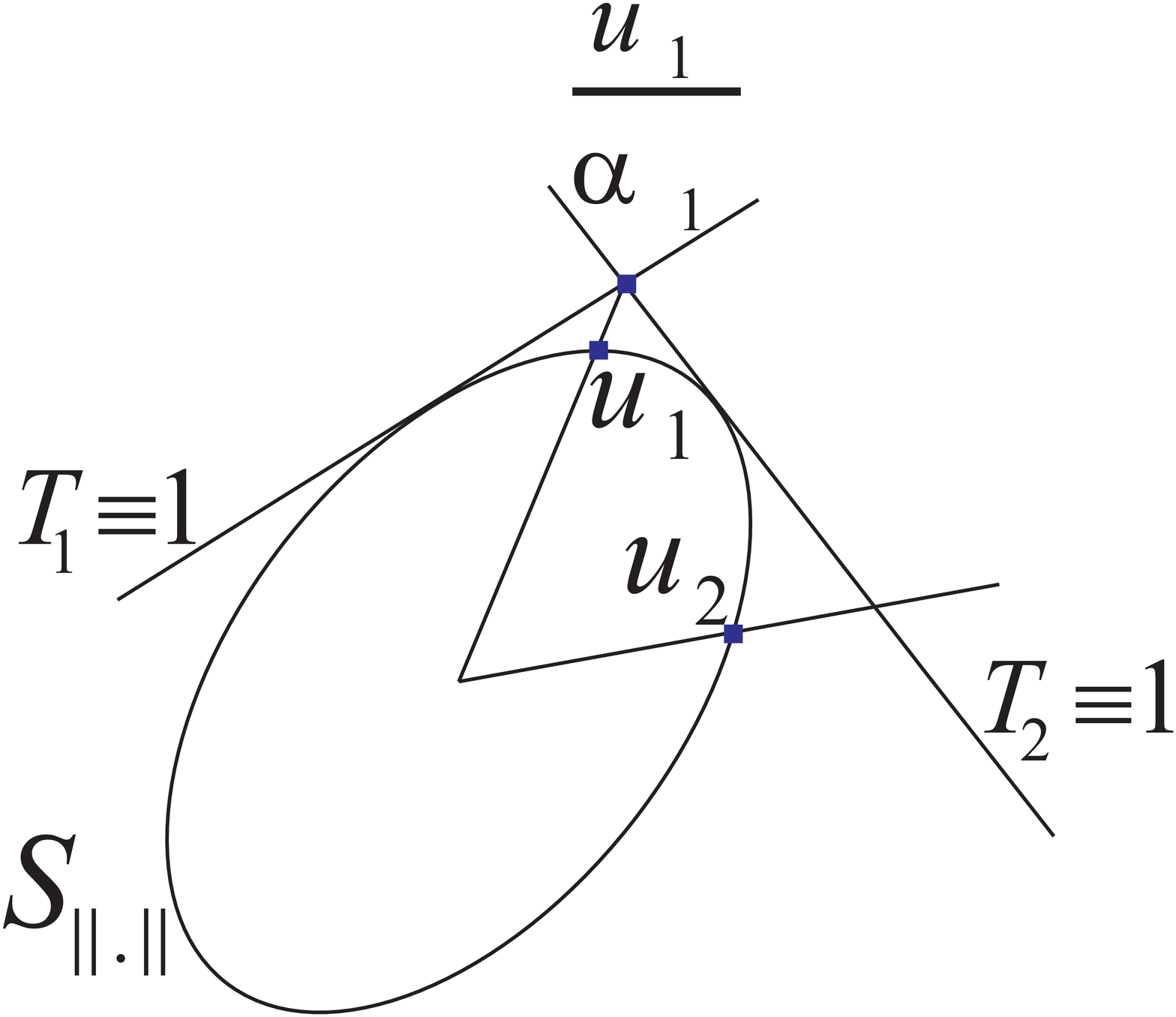} &
  \includegraphics[width=5cm]{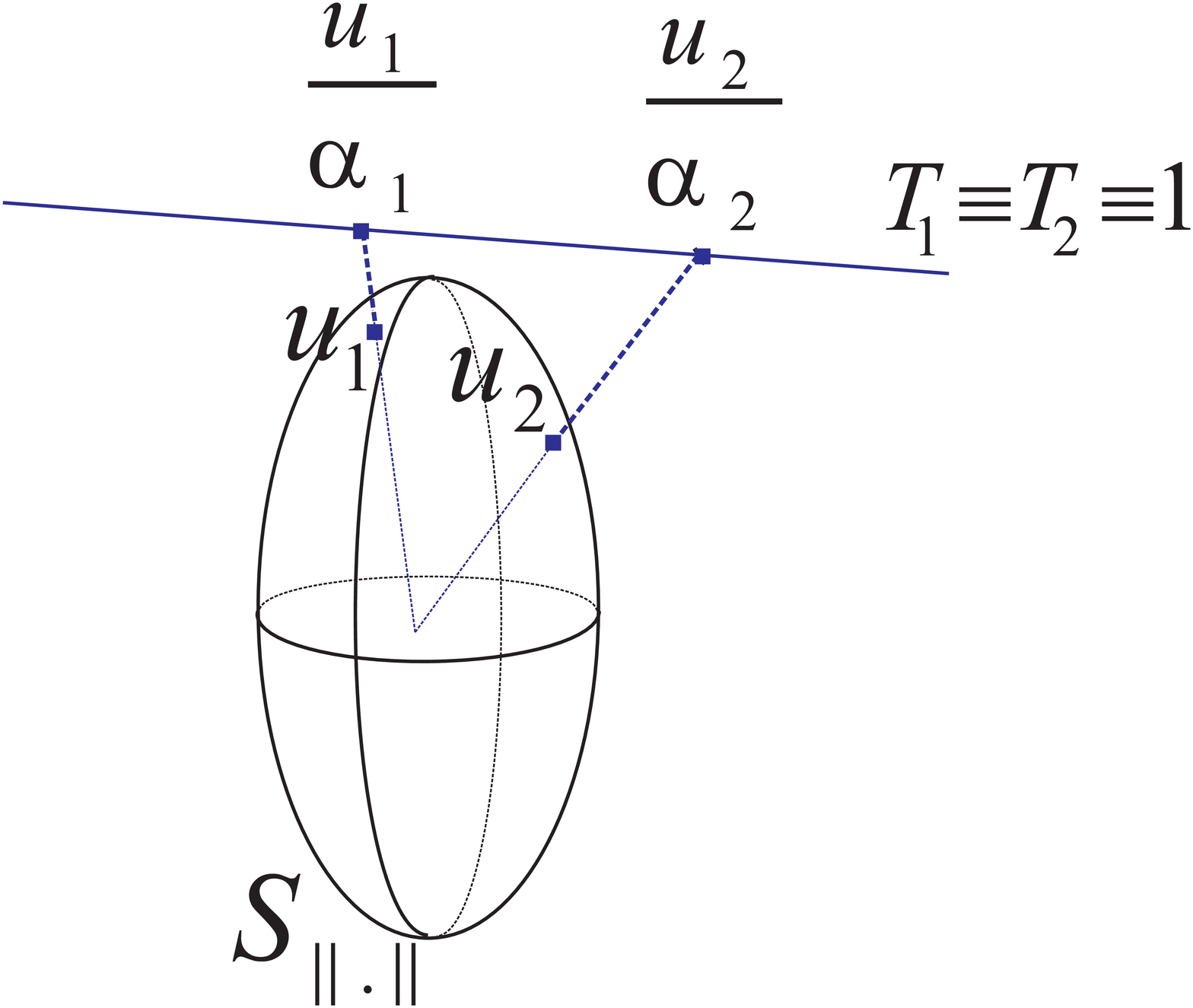} &
  \includegraphics[width=5cm]{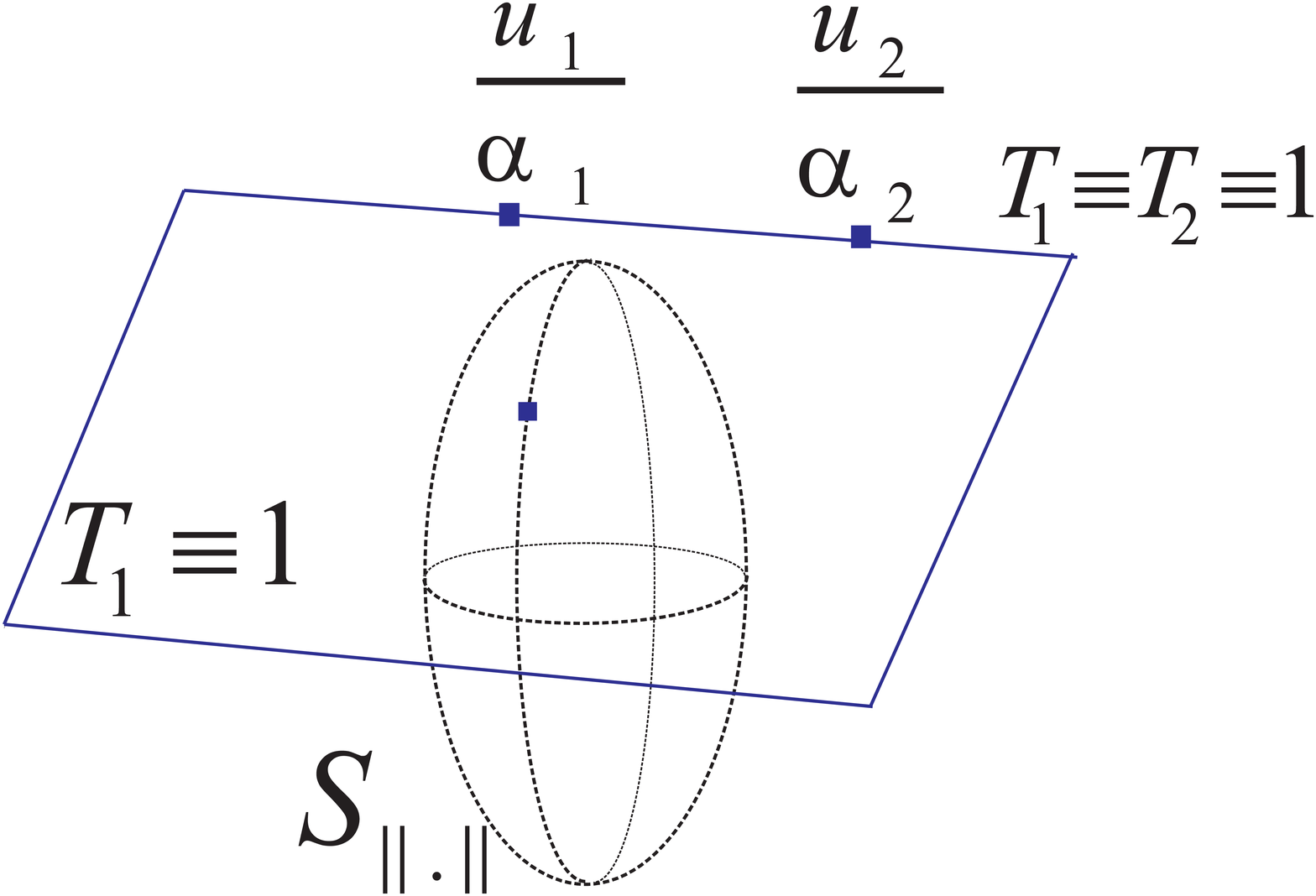} \\
  (a) & (b) & (c)
\end{tabular}
  \caption{(a): Case $n=2$. (b) and  (c):  Case $n=3$.}
  \label{}
\end{figure}

 As the above pictures for the
cases $n=2$ and $n=3$ suggest, there are at most two different
tangent affine hyperplanes to the unit sphere $S_{||\cdot||}$
containing the affine subspace passing through the points
$\frac{u_1}{\alpha_1}$,...,$\frac{u_{n-1}}{ \alpha_{n-1}}$, that
is there are at most two different linear mappings $T_1$, $T_2$ in
$F$. Indeed, assume first that the cardinal of $F$ is  one and let
us denote by $T_1$ the element of $F$. Then any real number
$\alpha \not\in\{0,\,T_1(u_n)\}$ satisfies condition \eqref{vacio}
and $\mathcal{R}=\mathbb R^*\setminus \{T_1(u_n)\}$. Now, if there
are two elements on $F$, $T_1\not=T_2$, we claim that any other
{\em different} element of $F$, say $T_3$, can be written as
$T_3=\gamma T_1+(1-\gamma)T_2$ for some $\gamma\in\mathbb
R\setminus\{0,1\}$. Indeed, since $T_3-T_2\not=0$ and
$T_1-T_2\not=0$ belong to the one dimensional subspace
$u_1^{\bot}\cap...\cap u_{n-1}^{\bot}\cap[u_1,...,u_n]^*$, we have
that there is $\gamma\in\mathbb R$ with $T_3-T_2=\gamma(T_1-T_2)$.
Since $T_3\not=T_2$, the constant $\gamma\not=0$. In addition,
since $T_3\not= T_1$, we have that $\gamma \not=1$, which proves
our claim.

Furthermore, this implies that the three {\em different} points
$T_1,\,T_2$ and $T_3$ of the dual unit sphere lie in a common
line, and, because $\|\cdot\|^{*}$ is convex, the segment line
which passes through $T_1,\,T_2$ and $T_3$, and whose end points
are two of these three points, is included in $S_{\|\cdot\|^*}$.
But this is in contradiction to the fact that the dual norm
$\|\cdot\|^*$ is strictly convex (because the norm $\|\cdot\|$ is
differentiable and $Z$ is finite dimensional, see \cite{DGZ}).
Finally, any real number $\alpha
\not\in\{0,\,T_1(u_n),\,T_2(u_n)\}$ satisfies condition
\eqref{vacio} and $\mathcal{R}=\mathbb
R^*\setminus\{T_1(u_n),\,T_2(u_n)\}$.
\end{proof}

\medskip

From the proof of Lemma \ref{case N} we obtain the following.

\begin{lem}\label{vectoradicional}
Let $Z=[u_1,...u_n]$ be a n-dimensional  space ($n\in \mathbb N $)
with a differentiable norm $||\cdot||$. Consider real numbers \
$0<\alpha_i<1$, \, for $i=1,...,n-1$. Then, the cardinal of the
set
\begin{equation}\label{vacio3x}
\big\{T\in S_{||\cdot||^*}: \ \ T(u_i)=\alpha_i\,,\
i=1,...,n-1\bigr\}.
\end{equation}
is at most two.
\end{lem}

\bigskip

The general strategy of the proof of Theorem \ref{approximation
theorem} is as follows. We consider the space \ $Y=X\oplus
\mathbb{R}$ \ and define the following norm on $Y$:
\begin{itemize}
\item  If \ $p>1$, \ for every $y=(x,r)\in Y$, put \
$N(y)=N(x,r)=(||x||^2+r^2)^{1/2}$, \ where $||\cdot||$ is a LUR
and \ $C^p$ \ smooth norm on \  $X$. Then, clearly the norm \ $N$
\ is
 \ LUR  \ and $C^1$ smooth on $Y$. Moreover, \ $N$ \ is \ $C^{p}$ \ smooth on \
 the open set $Y\setminus\{(0,\lambda):\ \lambda\in \mathbb R\}$.
Define $\nu=(0,1)$ and take  \ $\beta\in Y^*\setminus\{0\}$ \ such
that \ $X=\ker \beta$. \ Select \ $\beta_1\in Y^*\setminus[\beta]$
\ such that $\beta_1(\nu)\not=0$ \ and \  $\omega \in \ker
\beta\setminus \ker \beta_1 $. \ Consider the closed hyperplane of
$Y$, \ $X_1=\ker \beta_1 $. Then, the restriction of the norm $N$
to $X_1$ is a $C^p$ smooth and LUR norm  on $X_1$. Now, the
(equivalent) norm considered in \ $Y=X_1\oplus [\omega]$, \
defined as \ $|z+\lambda \omega|=(N(z)^2+\lambda^2)^{1/2}$, where
$z\in X_1$ \ and \ $\lambda\in \mathbb R$, \ is LUR  and $C^1$
smooth on $Y$ \ and $C^p$ \ smooth on \ $Y\setminus[\omega]$. In
particular, the norm \ $|\cdot|$ \ is \ $C^p$ \ smooth on the open
set \ $\mathcal{U}=Y\setminus \ker \beta=\{(x,r):\ x\in
X,\,r\not=0\}$.

It could also be proved that the Banach space $Y$ admits an
equivalent LUR and  $C^p$ smooth norm on $Y$ with bounded
derivatives up to the order $p$. Nevertheless, a LUR and $C^1$
smooth norm on $Y$ \ and \ $C^p$ \ smooth on \ $\mathcal{U}$, \ is
sufficient to prove our result. Recall that if $X$ has a LUR and
$C^p$ smooth norm \ and $p>1$, then $X$ is superreflexive
\cite{DGZ}.

\medskip

\item If $p=1$, since the dual space $Y^*$ is separable, there is
a norm $|\cdot|$ on $Y$ which is LUR, $C^1$ smooth and WUR whose
dual is strictly convex \cite{DGZ}. Recall that if the norm
$|\cdot|$ is WUR, then the dual norm $|\cdot|^*$ is uniformly
G\^{a}teux smooth, and thus, G\^{a}teaux smooth.
\end{itemize}

\medskip

 \noindent Therefore, if $X$ is reflexive, the dual norm
$|\cdot|^*$ is LUR and $C^1$ smooth  \cite{DGZ}. If $X$ is not
reflexive, the dual norm  $|\cdot|^*$ is strictly convex and
G\^{a}teaux smooth.

\medskip

Let us denote\ $S:=S_{|\cdot|}$, the unit sphere of $(Y,|\cdot|)$
\ and \ $S^*:=S_{|\cdot|^*}$, the unit sphere of
$(Y^*,|\cdot|^*)$. Let us consider, the duality mapping of the
norm $|\cdot|$ defined as
\begin{align*}
 D&: S \longrightarrow S^* \\ D&(x)=|\cdot|'(x), \end{align*} which   is $|\cdot|-|\cdot|^*$
continuous because the norm $|\cdot|$ is of class $C^1$.

\medskip

 We establish a $C^p$ diffeomorphism  $\Phi$ between $X$ and half
unit sphere in $Y$,\   $S^+:=\{y=(x,r)\in Y:\ r>0\}$,\ as follows:
$\Phi:X\longrightarrow S^+ $ \ is the composition \ $\Phi=\Pi\circ
i$, \ where $i$ is the inclusion
 \ $i:X\longrightarrow Y,\ i(x)=(x,1)$ \ and $\Pi$ is defined by \
 $\Pi:Y\setminus\{0\}\longrightarrow S$, \
$\Pi(y)=\frac{\,y\,}{\,|y|\,}$.

\medskip

In order to simplify the notation, we will make the proof for the
case of a constant $\varepsilon>0$. By taking some standard
technical precautions the same proof will work in the case of a
positive continuous function $\varepsilon:X\to (0,+\infty)$ (at
the end of the proof we will explain what small changes should be
made).

Now, given a continuous function \ $f:\,X
\longrightarrow\mathbb{R}$, \ we consider the composition \
$F:=f\circ\Phi^{-1}:\,S^+\longrightarrow \mathbb{R}$, which is
continuous as well.\ For any given $\varepsilon >0$ we will \
$3\varepsilon$-approximate \ $F$ \ by a $C^p$ smooth function \
$H:S^+\longrightarrow \mathbb{R}$ \ with the properties that:
\begin{itemize} \item the set of critical points of $H$ is the countable
union of a family of disjoint sets $\{K_n\}_n$\,, \item there are
countable families of open slices $\{O_n\}_n$ \ and \ $\{B_n\}_n$
in $S^+$, such that \ $\cup_n\overline{B}_n$ is relatively  closed
in $S^+$, \ $\dist(B_n,\ X\times\{0\})>0$, \ $K_n\subset
O_n\subset B_n$\,, \ $\dist(O_n,\,S^+\setminus B_n)>0$ \  and \
$\dist(B_n, \cup_{m\not=n}B_m)>0$, \ for every $n\in\mathbb N$,
\item the oscillation of $F$ in every $B_n$ is less than
$\varepsilon$.\end{itemize} \noindent (We will consider slices of
$S^+$ of the form $\{x\in S: f(x)>r\}$\,,\ where $f$ is a
continuous linear functional of norm one and $0<r<1$. Recall also
that the distance between two sets $A$ and $A'$ in a Banach space
$(M,|\cdot|_M)$ is defined as the real number \
$\dist(A,A'):=\inf\{|a-a'|_M:\ a\in A\,,\ a'\in A'\}$.)

\medskip

\noindent Then we will prove that the function $h:=H\circ\Phi$ is
a $C^p$ smooth function on $X$, which $3\varepsilon$-approximates
$f$, and the set of critical points of $h$, $C=\{x\in X:
h'(x)=0\}$, can be written as
$C=\bigcup_{n=1}^{\infty}\mathcal{K}_{n}$, where, for every $n\in
\mathbb N$, the set $\mathcal{K}_n:=\Phi^{-1}(K_n)$ is contained
in the {\em open, convex, bounded and $C^p$ smooth body}
$\mathcal{O}_n:=\Phi^{-1}(O_n)$, which in turn is contained in the
open, convex, bounded and $C^p$ smooth body
$\mathcal{B}_n:=\Phi^{-1}(B_n)$, in such a way that
$\textrm{dist}(\mathcal{O}_n, X\setminus \mathcal{B}_n)>0$, \ the
oscillation of $f$ in $\mathcal{B}_n$ is less than \
$\varepsilon$, \ $\cup_n\overline{\mathcal{B}}_n$ is closed \ and
\ $\dist(\mathcal{B}_n, \cup_{m\not=n}\mathcal{B}_m)>0$. Once we
have done this, we will compose the function $h$ with a sequence
of deleting diffeomorfisms which will eliminate the critical
points of $h$. More precisely, for each set
 $\mathcal{O}_n$ we will find a $C^p$ diffeomorphism $\Psi_n$ from $X$ onto
$X\setminus \overline{\mathcal{O}}_n$ so that $\Psi_n$ is the
identity off $\mathcal{B}_n$. Then, by defining
$g:=h\circ\bigcirc_{n=1}^{\infty}\Psi_{n}$, we will get a $C^p$
smooth function which $4\varepsilon$- approximates $f$ and which
has no critical points.

\medskip

The most difficult part in this scheme is the construction of the
function $H$. We will inductively define linearly independent
functionals $g_{k}\in Y^{*}$, open subsets $U_k$ of $S^{+}$,
points $x_k\in U_k$, real numbers $a_k\not=0$ satisfying
$|a_k-F(x_k)|<\varepsilon$, real numbers $\gamma_{k}$ and
$\gamma_{i,j}$  in the interval $(0,1)$ \ (with $i+j=k$),
functions $h_k$ of the form
$h_k=\varphi_k(g_k)\,\phi_{k-1,1}(g_{k-1})\,\cdots\phi_{1,k-1}(g_1),$
\ where the \ $\varphi_k$,\ $\phi_{k-1,1}$,\,...,\,$\phi_{1,k-1}$
are suitably chosen $C^\infty$ functions on the real line, and
functions $r_k$ of the form $r_k=s_k g_k+(1-s_kg_k(x_k))$ (with
very small $s_k\not=0$), and put
    $$
    {\bf H_k}=\frac{\sum_{i=1}^k a_ir_ih_i}{\sum_{i=1}^k h_i},
    $$
where ${\bf H_k}: U_1\cup...\cup U_k\longrightarrow \mathbb R$.
The interior of the support of $h_k$ will be the set
\begin{equation*}
U_k=\{x\in S^{+}:\
 g_1(x)<\gamma_{1,k-1},..., \, g_{k-1}(x)<\gamma_{k-1,1} \  \text{ and }\
 g_k(x)>\gamma_k\},
 \end{equation*}
where the oscillation of the function $F$ will be less than
$\varepsilon$. Denote by $T_x$ the (vectorial) tangent hyperplane
to $S^+$ at the point $x$, that is $T_x:=\operatorname{ker}D(x)$.
The derivative of ${\bf H_k}$ at every point $x\in U_1\cup...\cup
U_k$ will be shown to be the restriction to $T_x$ of a {\em
nontrivial linear combination of the linear functionals } $g_{1},
..., g_{k}$. Then, by making use of Lemmas \ref{case N} and
\ref{vectoradicional} and choosing the $\gamma_{i,j}$ close enough
to $\gamma_{i}$, we will prove that the set of critical points of
${\bf H_k}$ is a finite union of pairwise disjoint sets which are
contained in a finite union of pairwise disjoint slices, with
positive distance between any two slices (see Figure \ref{dibujo
n=3} below). These slices will be determined by functionals in
finite sets $N_k\subset Y^*$ defined by a repeated application of
Lemmas \ref{case N} and \ref{vectoradicional}. The function $H$
will be then defined as
    $$
    H=\frac{\sum_{k=1}^\infty a_kr_kh_k}{\sum_{k=1}^\infty h_k}.
    $$

\medskip

Let us begin with the formal construction of the functions ${\bf
H_{k}}$. We will use the notation $H_k$ and $H'_{k}$ when $H_k$
and its derivative $H_{k}'$ are thought to be defined on an open
subset of $Y$ and reserve the symbols ${\bf H_{k}}$ and ${\bf
H_{k}'}(x)$ for the restriction of $H_k$ and $H_{k}'(x)$ to a
subset of $S$ and to the tangent space $T_{x}$ of $S$ at $x$,
respectively.

\medskip

 Since the norm  $|\cdot|$ is LUR
 we can find, for every $x\in S^+$,  open slices $R_x=\{y\in S: \ f_x(y)>\delta_x\}\subset
 S^+$ \ and \  $P_x=\{y\in S: \ f_x(y)>\delta_x^4\}\subset S^+$,
 \ where \ $0<\delta_x<1$ \ and \ $|f_x|=1=f_x(x)$, \
 so that the oscillation of $F$ in every \ $P_x$ is less than
 $\varepsilon$. We also assume, for technical reasons, and with no
 loss of generality, that
 $\dist(P_x,\,X\times\{0\}\,)>0$.

\medskip

\noindent Since $Y$ is separable we  can select a countable
subfamily of $\{R_x\}_{x\in S^+}$, which covers $S^+$.  Let us
denote this countable subfamily by $\{S_n\}_n$, where $S_n:=\{y\in
S: f_n(y)>\delta_n\}$. Recall that the oscillation of $F$ in every
$P_n:=\{y\in S: f_n(y)>\delta_n^4\}$ is less than \ $\varepsilon$
\ and \ $\dist(P_n, \ X\times\{0\})>0$.

\medskip

\noindent $\bold{\diamond}$ {\bf For $\bold{k=1}$}, define
\begin{align*} h_1&:\,S^+\longrightarrow \mathbb R \\ h_1&=\varphi_1(f_1), \end{align*}
 where $\varphi_1$ is a
$C^\infty$ function on $\mathbb R$ satisfying
\begin{align*}
\varphi_1(t)&=0 \ \text{ if }\  t\le\delta_1 \\
\varphi_1(1)&=1 \notag \\
\varphi_1'(t)&>0 \  \text{ if } \ t>\delta_1\, . \notag
\end{align*}
Notice that the interior of the support of $h_1$ is the open set
$S_1$. Denote by $x_1$ the point of $S^+$ satisfying $f_1(x_1)=1$.
Now select $a_1\in \mathbb R^*=\mathbb R\setminus\{0\}$\ with \
$|a_1-F(x_1)|<\varepsilon$ and define the auxiliary function
\begin{align*}
&r_1:S^+\longrightarrow \mathbb R,\\
&r_1=s_1f_1+(1-s_1f_1(x_1)),\notag \end{align*} where we have
selected $s_1$ so that $a_1s_1>0$ and $|s_1|$ small enough so that
the oscillation of $r_1$ on $S_1$ is less than \
$\frac{\varepsilon}{\,|a_1|\,}$. \,  Notice that $r_1(x_1)=1$.
Define
\begin{align*}
&{\bf H_1}:S_1\longrightarrow \mathbb R\\
&{\bf H_1}=\frac{a_1r_1h_1}{h_1}=a_1r_1.\notag
\end{align*}
The function $\bold{H_1}$ is $C^p$ smooth in $S_1$ and  the set of
critical points of ${\bf H_1}$,
\begin{equation*}Z_1= \{x\in S_1:\,H'_1(x)=0 \ \text{ on
} \ T_x\}\end{equation*} consists of the unique point $x_1$.
Indeed, $\bold{H'_1}(x)=H'_1(x)|_{T_x}=a_1s_1f_1|_{T_x}\equiv 0$ \
iff $D(x)=f_1$. This implies that  \ $Z_1=\{x_1\}$. Now select
real numbers $\gamma_{1,1}'$,\ $t_1$ and $l_1$ \ such that \
$\delta_1<\gamma_{1,1}'<t_1<\l_1<1$ and define the open slices
\begin{equation*}
 O_{f_1}=\{x\in S: \ f_1(x)>l_1\} \quad \text{ and } \quad
B_{f_1}=\{x\in S: \ f_1(x)>t_1\}.
 \end{equation*}
Clearly the above sets satisfy that  \ $Z_1\subset O_{f_1}\subset
B_{f_1}\subset S_1$, \  $\dist(O_{f_1},\, S\setminus B_{f_1})>0$ \
and \ $\dist(B_{f_1},\, \{x\in S:\ f_1(x)\le \gamma_{1,1}'\})>0$.

 \medskip

\noindent In order to simplify the notation in the rest of the
proof, let us denote by \ $\gamma_1=\delta_1$, \ $U_1=R_1=S_1$, \
$g_1=f_1$,  \ $z_1=x_1$ \ and $\Gamma_1=N_1=\{g_1\}$. Let us
define $\sigma_{1,1}=a_1s_1$ and write ${\bf
H_1'}=\sigma_{1,1}\bold{g_1}$ on $U_1$ where $\bold{g_1}$ is the
restriction of $g_1$ to $T_x$ whenever we evaluate
$\bold{H_1'}(x)$. In addition, if $A\subset S$, we denote by
$A^c=S\setminus A$.

\bigskip

\noindent $\bold{\diamond}$ {\bf  For $\bold{k=2}$}. Let us denote
by $y_2\in S^{+}$ the point satisfying $f_2(y_2)=1$. If either
$\{g_1,\ D(y_2)=f_2\}$ are lineally dependent (this only occurs
when $g_1=f_2)$ \ or \ $g_1(y_2)=\gamma_1$, we use the density of
the norm attaining functionals (Bishop-Phelps Theorem) and the
continuity of $D$ to  modify $y_2$ and find $z_2\in S^{+}$ so that
$\{g_1,\ D(z_2):=g_2\}$ are l.i., \ $g_1(z_2)\not=\gamma_1$ \ and
\begin{equation*}
\{x\in S:f_2(x)>\delta_2^2\}\subset\{x\in S:
g_2(x)>{\nu_2}\}\subset\{x\in S:f_2(x)>\delta_2^3\},
\end{equation*}
for some ${\nu_2}\in(0,1)$. If \ $g_1(y_2)\not=\gamma_1$ \ and
$\{g_1,f_2\}$ are l.i., define $g_2=f_2$ and $z_2=y_2$. \ Then,
apply Lemma \ref{case N} to the 2-dimensional space $[g_1,g_2]$
with the norm $|\cdot|^*$ (the restriction to $[g_1,g_2]$ of the
dual norm $|\cdot|^* $ considered in $Y^*$) and the real number
$\gamma_1\in(0,1)$ to obtain ${\gamma_2}\in(0,1)$ close enough to
$\nu_2$ so that
\begin{equation*}
S_2=\{x\in S:f_2(x)>\delta_2\}\subset\{x\in S:
g_2(x)>{\gamma_2}\}\subset\{x\in S: f_2(x)>\delta_2^4\}=P_2
 \end{equation*}
 and
\begin{equation}\label{empty2}
\{T\in [g_1,g_2]^* :\   |T|=1, \ T(g_1)=\gamma_1 \ \text{ and } \
T(g_2)={\gamma_2} \}=\emptyset
\end{equation}
Recall that the norm $|\cdot|^*$ is G\^{a}teaux differentiable on
$Y^*$ and therefore the restriction of this norm to $[g_1,g_2]$,
which we shall denote by $|\cdot|^*$ as well, is a differentiable
norm on the space $[g_1,g_2]$ (G\^{a}teaux and Fr\'{e}chet notions
of differentiability are equivalent in the case of {\em convex}
functions defined on {\em finite-dimensional} spaces). Therefore,
we can apply Lemma \ref{case N} to the norm $|\cdot|^*$ in the
space $[g_1,g_2]$. In fact, the same argument works for any finite
dimensional subspace of $Y^*$ and we will apply Lemma \ref{case N}
in the next steps to  larger finite dimensional subspaces of
$Y^*$. Define the sets
\begin{align*}R_2&=\{x\in S: g_2(x)>\gamma_2\},  \ \ \text{ and }\\
U_2'=\{&x\in S:\ g_1(x)<\gamma_{1,1}'\,,\
g_2(x)>\gamma_2\}.\end{align*} Assume that $U_1\cap
R_2\not=\emptyset$ \ and consider the  set
$$M_2=D^{-1}([g_1,g_2])\cap U_2'\cap U_1.$$
In the case that $M_2=\emptyset$, we  select as $\gamma_{1,1}$ any
point in $(\gamma_1,\gamma_{1,1}')$. In the case that
$M_2\not=\emptyset$ and $\gamma_1<\inf\{g_1(x):\ x\in M_2\}$, we
select $\gamma_{1,1}$ \ so that \
$$\gamma_1<\gamma_{1,1}<\inf\{g_1(x):\ x\in M_2\}.$$

\medskip

\noindent In the case that $\gamma_1=\inf\{g_1(x):\ x\in M_2\}$
and in order to obtain an appropriate $\gamma_{1,1}$ we need to
study the limits of the sequences $\{x_n\}\subset M_2$ such that
$\lim_n g_1(x_n)=\gamma_1$.  Define
\begin{equation*}
F_2'=\{T\in[g_1,g_2]^*: \ |T|=1  \ \text{ and } \ T(g_1)=\gamma_1
\}.\end{equation*} From Lemma \ref{vectoradicional}, we deduce
that the cardinal of the set $F_2'$ is at most two. Furthermore,
since $|\cdot|^*$ is strictly convex, the cardinal of  the set
\begin{equation*}
N_2'=\{g\in S^*\cap[g_1,g_2]: \ T(g)=1 \ \text{ for some } \ T\in
F_2'\}\end{equation*} is at most two.

\medskip

Let us take any sequence $\{x_n\}\subset M_2$ with
$\lim_ng_1(x_n)=\gamma_1$. Consider every $x_n$ as an element of
$X^{**}$ and denote by $\bold{x_n}$  its restriction to
$[g_1,g_2]$. Recall that if $x_n\in M_2$, then  $D(x_n)\in S^*
\cap[g_1,g_2]$, \ for every $n\in \mathbb N$.\  Moreover, the
sequence of restrictions $\{\bold{x_n}\}\subset [g_1,g_2]^*$
satisfies that
\begin{equation*}
1=|x_n|\ge|\bold{x_n}|=\max\{\bold{x_n}(h):\ h\in
S^*\cap[g_1,g_2]\}\ge
\bold{x_n}(D(x_n))=D(x_n)(x_n)=1,\end{equation*} for every $n\in
\mathbb{N}$. Thus, there is a subsequence \ $\{\bold{x_{n_j}}\}$ \
converging to an element \ $T\in [g_1,g_2]^*$ \ with \ $|T|=1$.
Since $\lim_jg_1(x_{n_j})=\lim_{j}\bold{x_{n_j}}(g_1)=\gamma_1$,
then $T(g_1)=\gamma_1$ and this implies that $T\in F_2'$.
Furthermore, if $g\in N_2'$ and $T(g)=1$, then
$\lim_j\bold{x_{n_j}}(g)=1$. In addition,
$T(g_2)=\lim_j\bold{x_{n_j}}(g_2)=\lim_j g_2(x_{n_j})\ge
\gamma_2$. Then, from condition \eqref{empty2}, we deduce that
$T(g_2)>\gamma_2$. Let us define
\begin{align*}
F_2=\{T\in F_2': \text{ there} & \text{ is a sequence }  \
\{x_n\}\subset M_2 \ \text{ with } \ \lim_n\bold{x_n}=T \ \text{
and } \ \lim_n\bold{x_n}(g_1)=\gamma_1 \},\\ N_2&=\{g\in N_2': \
T(g)=1 \ \text{ for some } \ T\in F_2\}.
\end{align*}
Select  a real number \ $\gamma_{2}'$ \  satisfying  \ $\gamma_2<
\gamma_{2}'<\min\{T(g_2):\ T\in F_2\}$ (recall that $F_2$ is
finite). Let us prove the following Fact.

\medskip

\begin{fact}\label{fact 1}
\hfill{ }
\begin{enumerate} \item There are numbers \ $0<t_2<l_2<1$ \ such that for every $g\in N_2$,
the slices $$O_g:=\{x\in S:\ g(x)>l_2\} \  \text{ and } \quad
B_g:=\{x\in S:\ g(x)>t_2\}$$ satisfy that
\begin{align}\label{inclusionfact1} O_g&\subset B_g\subset \{x\in S: \
g_1(x)<\gamma_{1,1}'\,,\ g_2(x)>\gamma_{2}'\} \ \text{ and }\\
\label{intersectionfact1} &\dist(B_g,B_{g'})>0, \text{ whenever }
g,g'\in N_2\,, \ g\not=g'.
\end{align}
\item There is  \  $\gamma_{1,1}\in (\gamma_1,\gamma_{1,1}')$,\
such that if $x\in M_2$ and $g_1(x)<\gamma_{1,1}$,\ then $x\in
O_g$, for some $g\in N_2$.
\end{enumerate}
\end{fact}
\noindent {\em Proof of Fact \ref{fact 1}.} (1) First, if $X$ is
reflexive, we know that for every \ $g\in N_2$ \  there is \
$x_g\in S$ \ such that \ $D(x_g)=g$.\  Since \ $\mathbf{x_g}(g)=1$
\ and \ $|\cdot|^*$ \ is G\^{a}teaux smooth, \ then $\bold{x_g}\in
F_2$. This implies that $\bold{x_g}(g_1)=\gamma_1<\gamma_{1,1}'$
and $\bold{x_g}(g_2)>\gamma_{2}'$. Hence, $x_g\in \{x\in S: \
g_1(x)<\gamma_{1,1}',\ g_2(x)>\gamma_{2}'\} $. Now, since  the
norm $|\cdot|$ is LUR and $D(x_g)=g$, the functional $g$ strongly
exposes $S$ at the point $x_g$. Taking into account that $N_2$ is
finite we can hence obtain real numbers $0<t_2<l_2<1$ and slices
$O_g$ and $B_g$ satisfying conditions \eqref{inclusionfact1} and
\eqref{intersectionfact1}, \ for every $g\in N_2$.

\medskip

\noindent Now consider a non reflexive Banach space $X$. Let us
first prove \eqref{inclusionfact1}. Assume, on the contrary, that
there is a point $g \in N_2$ and there is a sequence
$\{y_n\}\subset S$ satisfying $g(y_n)>1-\frac1n$ with either \
$g_1(y_n)\ge \gamma_{1,1}' $ \ or \ $g_2(y_n)\le \gamma_{2}'$ \
for every $n\in \mathbb N$. Since \ $g\in N_2$ \ there is a
sequence \ $\{x_n\}\subset M_2$ \ with \ $\lim_n
g_1(x_n)=\gamma_1$,\ $\lim_n g_2(x_n)>\gamma_{2}'$ \ and \ $\lim_n
g(x_n)=1$. In particular,
\begin{equation*}
\frac{g(x_n)+1-\frac1n}{2}< g\left(\frac{x_n+y_n}{2}\right)\le
\left|\frac{x_n+y_n}{2}\right|\le 1,
\end{equation*}
and thus $\lim_n\left|\frac{x_n+y_n}{2}\right|= 1$. Recall that,
in this case the norm $|\cdot|$ is WUR, and hence
$x_n-y_n\xrightarrow{\omega} 0$\  (weaky converges to zero). This
last assertion gives a contradiction since either \
$\limsup_ng_1(x_n-y_n)\le \gamma_1-\gamma_{1,1}'<0$ \ or \
$\liminf_ng_2(x_n-y_n)\ge \lim_ng_2(x_n)-\gamma_{2}'>0$. \
Therefore we can find real numbers \ $0<t_2<l_2<1$ \ and  slices \
$O_g$ \ and \ $B_g$ \ for every \ $g\in N_2$, satisfying condition
\eqref{inclusionfact1}. In order to obtain
\eqref{intersectionfact1} we just need to modify  $t_2$ and $l_2$
and select them close enough to $1$. Indeed, assume on the
contrary,  that there are sequences $\{y_n\}\subset S$ and
$\{z_n\}\subset S$ and $g,g'\in N_2$, $g\not=g'$, such that
$\lim_ng(y_n)=1$, \ $\lim_ng'(z_n)=1$ and $\lim_n|y_n-z_n|=0$.
Then,
\begin{align*}\frac{g(y_n)+g'(z_n)}{2}&\le
\frac{(g+g')(y_n)+g'(z_n-y_n)}{2}\le
\frac{(g+g')(y_n)+|z_n-y_n|}{2}\\&\le \frac{|g+g'|^*}{2}+
\frac{|z_n-y_n|}{2}\le 1+\frac{|z_n-y_n|}{2}.
\end{align*}
Since the limit of the first and last terms in the above chain of
inequalities is $1$, we deduce that \ $|g+g'|^*=2$. Since the norm
$|\cdot|^*$ is strictly convex, we deduce that $g=g'$, a
contradiction.

\medskip

 \noindent (2)  Assume, on the contrary, that for every $n\in \mathbb N$, there is \ $x_n\in
M_2$ \ with \ $g_1(x_n)\le \gamma_1+\frac1n$ \ and \ $\{x_n: n \in
\mathbb N\} \cap (\cup_{g\in N_2}O_g)=\emptyset$. Since \
$\lim_ng_1(x_n)=\gamma_1$ \ and \ $\{x_n\}\subset M_2$, \ from the
comments preceding Fact \ref{fact 1}, we know that there is a
subsequence $\{x_{n_j}\}$ and $g\in N_2$ satisfying that
$\lim_jg(x_{n_j})=1$, which is a contradiction. This finishes the
proof of Fact \ref{fact 1}. $\Box$

\medskip

\medskip

If  $U_1\cap R_2=\emptyset$, we  can select as $\gamma_{1,1}$ any
number in $(\gamma_1,\gamma_{1,1}')$. Now, we define,
\begin{align*}h_2&:\,S^+\longrightarrow \mathbb R\\
h_2&=\varphi_2(g_2)\,\phi_{1,1}(g_1)
\end{align*}
where $\varphi_2$ and $\phi_{1,1}$ are $C^\infty$ functions on
$\mathbb{R}$ satisfying:
\begin{align*}
\varphi_2(t)&=0 \ \text{ if }\  t\le\gamma_2 \\
\varphi_2(1)&=1 \\
\varphi_2'(t)&>0 \  \text{ if } \ t>\gamma_2\, ,
\end{align*}
and
\begin{align*}
\phi_{1,1}(t)&=1 \ \ \text{ if }\ t\le \textstyle{\frac{\gamma_1+\gamma_{1,1}}{2}} \\
\phi_{1,1}(t)&=0 \ \  \text{ if }\ t\ge \gamma_{1,1} \notag\\
\phi'_{1,1}(t)&<0 \ \  \text{ if }\ t\in
(\,\textstyle{\frac{\gamma_1+\gamma_{1,1}}{2}},\, \gamma_{1,1}).
\notag\end{align*}  Notice that the interior of the support of
$h_2$ is the open set $$U_2=\{x\in S: g_1(x)< \gamma_{1,1}\,,\
g_2(x)>{\gamma_{2}}\}.$$

\noindent Select one point $x_2\in U_2$, a real number $a_2\in
\mathbb R^*$ with $|a_2-F(x_2)|<\varepsilon$ and define the
auxiliary function
\begin{align*}
&r_2:S^+\longrightarrow \mathbb R,\\
&r_2=s_2g_2+(1-s_2g_2(x_2)),\notag \end{align*} where we have
selected $s_2$ so that $s_2a_2>0$ and $|s_2|$ is small enough so
that the oscillation of $r_2$ on $U_2$ is less than \
$\frac{\varepsilon}{|a_2|}$. \, Notice that $r_2(x_2)=1$.

\noindent Let us study the critical points $Z_2$ of the function
\begin{align*} &{\bf H_2}: U_1\cup U_2\longrightarrow \mathbb R,\\ \notag
&{\bf H_2}=\frac{a_1r_1h_1+a_2r_2h_2}{h_1+h_2}.
\end{align*}
Let us prove that $Z_2=\{x\in U_1\cup U_2:\,H'_2(x)= 0 \ \text{ on
} \ T_x\}$  can be included in a   finite number of pairwise
disjoint slices within $U_1\cup U_2$ by splitting it conveniently
into up to four sets.

\noindent First, \ if \ $x\in U_1\setminus U_2$,\ we have that
$H_2(x)=a_1r_1(x)$ \ and \
$\bold{H_2'}(x)=H_2'(x)|_{T_x}=a_1s_1g_1|_{T_x}\equiv 0$ \ iff \
$D(x)=g_1$. \ Thus, \ $Z_2\cap(U_{1}\setminus
U_{2})\subseteq\{z_{1}\}$. \ Second, if $x\in U_2\setminus U_1$,\
we have $H_2(x)=a_2r_2(x)$ \ and \
$\bold{H_2'}(x)=H_2'(x)|_{T_x}=a_2s_2g_2|_{T_x}\equiv 0$ iff
$D(x)=g_2$. \ Then, if  $z_2\in U_2\setminus U_1$, \ ${\bf H_2}$
has one critical point in $U_2\setminus U_1$, namely $z_2$; in
this case, since $g_1(z_2)\not=\gamma_1$, the point $z_2$ actually
belongs to $U_2\setminus \overline{U}_1$.

\noindent  Now, let us study the critical points of $\bold{H_2}$
in $U_1\cap U_2$. In order to simplify the notation, let us put
$\Lambda_1=\frac{h_1}{h_1+h_2}$, and denote by \ $\bold{g_1}$ \
and $\bold{g_2}$ \ the restrictions \ $g_1|_{T_x}$ \ and \
$g_2|_{T_x}$, respectively, whenever we consider $\bold{H_2'}(x)$
\ and $\Lambda_1'(x)$. Then, \ ${\bf H_2}=a_1r_1\Lambda_1
+a_2r_2(1-\Lambda_1)$ \ and
\begin{align*}
{\bf H_2'}&=a_1s_1\Lambda_1
\bold{g_1}+a_2s_2(1-\Lambda_1)\bold{g_2}
+(a_1r_1-a_2r_2)\Lambda'_1 \\  &=\sigma_{1,1}\Lambda_1
\bold{g_1}+a_2s_2(1-\Lambda_1)\bold{g_2}
+(H_1-a_2r_2)\Lambda'_1\end{align*} By computing $\Lambda_1'$, we
obtain $\Lambda_1'=\xi_{1,1}\bold{g_1}+ \xi_{1,2}\bold{g_2}$,
where the coefficients $\xi_{1,1}$ \ $\xi_{1,2}$ are continuous
functions on $U_1\cup U_2$ and have the following form,
\begin{align*}
 \xi_{1,1}&=\frac{\varphi_1'(g_1)h_2-h_1\varphi_2(g_2)\phi_{1,1}'(g_1)}{(h_1+h_2)^2} \\
\xi_{1,2}&=\frac{-h_1\varphi_2'(g_2)\,\phi_{1,1}(g_1)}{(h_1+h_2)^2}
\end{align*}
Thus \  $ {\bf
H_2'}=\sigma_{2,1}\bold{g_1}+\sigma_{2,2}\bold{g_2}$, \ where \
$\sigma_{2,1}$ and $\sigma_{2,2}$ are continuous functions on
$U_1\cup U_2$ and have the following form
\begin{align}\label{derivada de H2}
\sigma_{2,1}&=\sigma_{1,1}\Lambda_1+(H_1-a_2r_2)\xi_{1,1}
\\
\sigma_{2,2}&= a_2s_2(1-\Lambda_1)+(H_1-a_2r_2)\xi_{1,2}. \notag
\end{align}

\noindent Notice that if $x\in U_1\cap U_2$\,,  then \
$\sigma_{1,1}>0$, \ $a_2s_2>0$, \ $\Lambda_1>0$, \
$1-\Lambda_1>0$, \ $\xi_{1,1}>0$ \ and \ $\xi_{1,2}<0$. Therefore,
on $U_1\cap U_2$\,, the coefficient $\sigma_{2,1}$ is strictly
positive whenever $H_1-a_2r_2\ge 0$, \ and the coefficient
$\sigma_{2,2}$ is strictly positive whenever $H_1-a_2r_2\le 0$.
Since the vectors $g_1$ and $g_2$ are l.i., if $x\in U_1\cap U_2$
\ and ${\bf H_2'}(x):T_x\longrightarrow \mathbb R$ \ is
identically zero, \ there is necessarily \ $\varrho\not=0$ with
$D(x)=\varrho (\sigma_{2,1}(x)g_1+\sigma_{2,2}(x)g_2)$. Thus,
$D(x)\in [g_1,g_2]$.

\medskip

\noindent The set $Z_2$ can be split into the disjoint sets
$Z_2=Z_1 \cup Z_{2,1}\cup Z_{2,2}$, where
\begin{equation*}
Z_{2,1}=\begin{cases}\{z_2\}, & \text{ if } \ z_2\in U_2\setminus
\overline{U}_1\\ \emptyset, & \text{ otherwise } \end{cases}
\end{equation*}
and \ $Z_{2,2}$ \ is a subset (possibly empty) within  \ $ U_1\cap
U_2 \cap D^{-1}([g_1,g_2])$. Now, let us check that
$Z_{2,2}\subset\cup_{g\in N_2}O_g$. \ Indeed, if \ $x\in
Z_{2,2}$\,, \ then \ $x\in U_1\cap U_2\subset U_1\cap U_2'$\,, \
$D(x)\in [g_1,g_2]$ \ and $g_1(x)<\gamma_{1,1}$. \  This implies,
according to Fact \ref{fact 1}, that \ $x\in \cup_{g\in N_2}O_g$.

\medskip

In the case when \ $Z_{2,1}=\{z_2\}$ \ and \ $z_2\not\in
\cup_{g\in N_2}{\overline{O}_g}$\,,  we select, if necessary, a
larger $t_2$\, \ with $t_2<l_2$\,,
 so that \ $z_2\notin \cup_{g\in N_2}\overline{B}_g$.
 Since the norm \  $|\cdot|$ \
is LUR and \ $D(z_2)=g_2$\,, \ the functional $g_2$ strongly
exposes $S$ at the point $z_2$ and we may select numbers \
$0<t_2'<l_2'<1$ \ and open slices, which are neighborhoods of \
$z_2$\,, \ defined by
\begin{equation*}
 O_{g_2}:=\{x\in S: g_2(x)>l_2'\} \quad \text{ and } \quad
B_{g_2}:=\{x\in S: g_2(x)>t_2'\},
\end{equation*}
  satisfying  \ $O_{g_2}\subset B_{g_2}\subset \{x\in S: \
g_1(x)<\gamma_{1,1}',\ g_2(x)>\gamma_{2}'\} $ and
$\dist(B_{g_2},B_g)>0$ for every $g\in N_2$. In this case, we
define $\Gamma_2=N_2\cup\{g_2\}$.

\medskip

 Now, if $Z_{2,1}=\{z_2\}\in \cup_{g\in N_2}{\overline{O}_g}$,\ we
 select, if necessary, a smaller constant $l_2$, with $0<t_2<l_2<1$,  so that
 $Z_{2,1}=\{z_2\}\in \cup_{g\in N_2}{O_g}$\,. In this case, and also when $Z_{2,1}=\emptyset$,
 we define $\Gamma_2=N_2$.

 Notice that, in any of the cases mentioned above, Fact \ref{fact 1} clearly
 holds for the (possibly) newly selected real numbers $t_2$ and $l_2$.

\medskip

 Notice that the distance between  any two sets  $B_{g}$, \  $B_{g'}$, \ $g,g'\in
 \Gamma_1\cup \Gamma_2$, \
 $g\not=g'$, \  is positive. \
Moreover, $Z_1\subset O_{g_1}\subset B_{g_1}\subset U_1=R_1$, \
and \ $Z_{2,1}\cup Z_{2,2}\subset \cup_{g\in \Gamma_2}O_g \subset
\cup_{g\in \Gamma_2} B_g \subset U_2'\subset R_2$. Therefore,
$Z_2=Z_1\cup Z_{2,1}\cup Z_{2,2}\subset \cup_{g\in \Gamma_1 \cup
\Gamma_2}O_g \subset \cup_{g\in \Gamma_1 \cup \Gamma_2}B_g \subset
U_1\cup U_2=R_1\cup R_2$. In addition, we have $\dist(\cup_{g\in
\Gamma_1 \cup \Gamma_2}B_g , \ (U_1\cup U_2)^c)>0$.

\medskip

It is worth remarking that ${\bf
H_2'}=\sigma_{2,1}\bold{g_1}+\sigma_{2,2}\bold{g_2}$ in $U_1\cup
U_2$, where $\sigma_{2,1}$ and $ \sigma_{2,2}$ are continuous
functions and at least one of the coefficients
$\sigma_{2,1}(x),\,\sigma_{2,2}(x)$ is strictly positive, for
every \ $x\in U_1\cup U_2$. Moreover, $\sigma_{2,1}(x)=0$ whenever
$x\not\in U_1$, and $\sigma_{2,2}(x)=0$ whenever $x\not\in U_2$.

\medskip

In order to clarify the construction in the general case, let us
also explain in detail the construction of the function $h_3$ and
locate the critical points of the function $\bold{H_3'}$.

\medskip

\noindent $\bold{\diamond}$ {\bf For $\bold{j=3}$,} let us denote
by $y_3\in S$ the point satisfying $f_3(y_3)=1$. If
$\{g_1,g_2,f_3\}$ are lineally dependent, \ or if \
$g_1(y_3)=\gamma_1$, \ or if \ $g_2(y_3)=\gamma_2$\,, we can use
the density of the norm attaining functionals (Bishop-Phelps
Theorem) and the continuity of $D$ to  modify $y_3$ and find
$z_3\in S$ so that: $g_1(z_3)\not=\gamma_1$,\
$g_2(z_3)\not={\gamma_2}$,\ $\{g_1,g_2,g_3:=D(z_3)\}$ are linearly
independent (l.i.), and
\begin{equation*}
\{x\in S:\ f_3(x)>\delta_3^2\}\subset\{x\in S:\
g_3(x)>{\nu_3}\}\subset\{x\in S: f_3(x)>\delta_3^3\}
\end{equation*}
for some ${\nu_3}\in (0,1)$. If $\{g_1,g_2,f_3\}$ are l.i., \
$g_1(y_3)\not=\gamma_1$, \ and   \ $g_2(y_3)\not=\gamma_2$, we
define $g_3=f_3$ and $z_3=y_3$. Then, we apply Lemma \ref{case N}
to the l.i. vectors, $\{g_1,g_2,g_3\}$ and the real numbers
$\gamma_1\in(0,1)$,\ ${\gamma_2}\in(0,1)$ and obtain $\gamma_3\in
(0,1)$ close enough to $\nu_3$ so that
\begin{equation*}
S_3=\{x\in S:\ f_3(x)>\delta_3 \}\subset\{x\in S:\
g_3(x)>{\gamma_3}\}\subset \{x\in S: \ f_3(x)>\delta_3^4\}=P_3,
\end{equation*}
\begin{equation}\label{3-1}
\{T\in [g_1,g_2,g_3]^*:\ T(g_1)=\gamma_1\,, \
T(g_2)={\gamma_2}\,,\ T(g_3)={\gamma_3}\ \text{ and } \ |T|=1
\}=\emptyset,
\end{equation}
\begin{equation}\label{3-2}
\{T\in [g_1,g_3]^*:\ T(g_1)=\gamma_1\,,\ T(g_3)={\gamma_3}\ \text{
and } \ |T|=1 \}=\emptyset,
\end{equation}
\begin{equation}\label{3-3}
\{T\in [g_2,g_3]^*:\ \ T(g_2)={\gamma_2}\,,\ T(g_3)={\gamma_3}\
\text{ and } \ |T|=1 \}=\emptyset.
\end{equation}

Select $\gamma_{2,1}'\in (\gamma_2,\gamma_{2}' )$ \ and define
\begin{align*} R_3=&\{x \in S:\, g_3(x)>\gamma_3\}  \ \text{ and } \\
U_3'=\{s\in  S: g_1&(x)<  \gamma_{1,2}'\,,\ g_2(x)<\gamma_{2,1}'\
\text{ and } \ g_3(x)>\gamma_3\},\end{align*} where
$\gamma_{1,2}'$ is a number in
$(\gamma_1,\,\frac{\gamma_1+\gamma_{1,1}}{2}).$

\noindent Notice that $\operatorname{dist}(B_g,\,U_3')>0$ \ for
every \ $g\in \Gamma_1\cup \Gamma_2$\,. Assume that $R_3\cap
(U_1\cup U_2)\not = \emptyset$, and consider the sets
\begin{align*}
M_{3,1}&=\{x\in (U_1\cap U_3') \setminus U_2:\  D(x)\in[g_1,g_3]\},\\
M_{3,2}&=\{x\in (U_2\cap U_3') \setminus U_1:\  D(x)\in[g_2,g_3]\},\\
M_{3,1,2}&=\{x\in U_1\cap U_2 \cap U_3':\ D(x)\in[g_1,g_2,g_3]
\},\end{align*} and $M_3=M_{3,1}\cup M_{3,2}\cup M_{3,1,2}$\,.

In the case that $M_3=\emptyset$, we select as $\gamma_{2,1}$ any
point in $(\gamma_2,\gamma_{2,1}')$ and $\gamma_{1,2}$ any point
in $(\gamma_1,\gamma_{1,2}' )$.

In the case that $M_3\not=\emptyset$ and $\dist(M_3, (U_1\cup
U_2)^c)>0$ we can easily find $\gamma_{2,1}\in
(\gamma_2,\gamma_{2,1}')$ and $\gamma_{1,2}\in
(\gamma_1,\gamma_{1,2}' )$ with $M_3\subset \{x\in S:
g_1(x)>\gamma_{1,2}\}\cup  \{x\in S: g_2(x)>\gamma_{2,1}\}$.

In the case that \  $\dist(M_3, (U_1\cup U_2)^c)=0$ \ and in order
to obtain suitable constants $\gamma_{2,1}$ and $\gamma_{1,2}$\,,
we need to study the limits of the sequences $\{x_n\}\subset M_3$
such that $\lim_n \dist(x_n,(U_1\cup U_2)^c)=0$. Define the sets
\begin{align*} F_{3,i}'&=\{T\in[g_i,g_3]^*:
\ T(g_i)=\gamma_i \ \text{ and } \ |T|=1\} \quad \text{ for } i=1,2, \\
 F_{3,1,2}'&=\{T\in[g_1,g_2,g_3]^* :\ T(g_1)=\gamma_1, \
T(g_2)=\gamma_2 \ \text{ and } \ |T|=1 \}, \end{align*} and
\begin{align*} N_{3,i}'&=\{g\in S^*\cap [g_i,g_3]:
\ T(g)=1 \  \text{ for some  } \ T\in F_{3,i}'\} \quad \text{ for } i=1,2,\\
N_{3,1,2}'&=\{g\in S^* \cap[g_1,g_2,g_3] :\ T(g)=1 \ \text{ for
some } \ T\in F_{3,1,2}'\}. \end{align*} Since the norm
$|\cdot|^*$ is
 G\^{a}teaux smooth, we apply Lemma \ref{vectoradicional} to the finite dimensional space
 $[g_1,g_2,g_3]$ and the restriction of the norm $|\cdot|^*$ to  $[g_1,g_2,g_3]$ (which is a differentiable norm on
 the space $[g_1,g_2,g_3]$), and deduce
that the cardinal of any of the sets $F_{3,i}'$\,, $F_{3,1,2}'$ is
at most two. Furthermore, from the strict convexity of the norm
$|\cdot|^*$ we obtain that the cardinal of any of the sets
$N_{3,i}'$ and $N_{3,1,2}'$\,, \ is at most two. Let us consider,
for $i=1,2$, the norm-one extensions to \ $[g_1,g_2,g_3]$ \ of the
functionals of $F_{3,i}'$\,, that is,
\begin{equation*} F_{3,i}''=\{T\in [g_1,g_2,g_3]^*:\, T|_{[g_i,g_3]}\in F_{3,i}' \ \text{ and } \ |T|=1\}. \end{equation*}
Since the norm $|\cdot|^*$ is G\^{a}teaux smooth, for every $G\in
F_{3,i}'$ there is exactly {\em one} norm-one extension $T$ to
$[g_1,g_2,g_3]$. Therefore the cardinal of the set $F_{3,i}''$ is
at most two. Hence the sets $F_3':=F_{3,1}''\cup F_{3,2}''\cup
F_{3,1,2}'$ and $N_3':=N_{3,1}'\cup N_{3,2}'\cup N_{3,1,2}'$  are
finite. In addition, as a consequence of the equalities
\eqref{3-1}, \eqref{3-2} and \eqref{3-3}, we deduce that
$T(g_3)\not=\gamma_3$ \ for every \ $T\in F_3'$. Indeed,  if $T\in
F_{3,1,2}'$ the assertion follows immediately   from  \eqref{3-1}.
If \ $T\in F_{3,i}''$ \ for some \ $i\in\{1,2\}$, \ then \
$T|_{[g_i,g_3]}\in F_{3,i}'$, that is, $|T|_{[g_i,g_3]}|=1$ \ and
\ $T(g_i)=\gamma_i$. From \eqref{3-2} \ for $i=1$, and \eqref{3-3}
for $i=2$, we obtain that $T(g_3)\not=\gamma_3$. We can restrict
our study to one of the following kind of sequences:

\medskip

\begin{enumerate}
\item Fix $i\in \{1,2\}$. Consider any sequence $\{x_n\}\subset
M_{3,i}$ \ such that \ $\lim_n \dist(x_n,(U_1\cup U_2)^c)=0$.
Then, it easily follows that $\lim_n g_i(x_n)=\gamma_i$. Indeed,

\begin{itemize}
\item if $\{x_n\}\subset M_{3,1}$, then in particular \
$\{x_n\}\subset U_1=R_1$. Therefore, \ $\dist(x_{n}, (U_1\cup
U_2)^c)\ge \dist(x_{n}, R_1^c)$. Thus, $\lim_n \dist(x_{n}, R_1^c)
=0$ and this implies that $\lim_n g_1(x_n)=\gamma_1$;

\item if $\{x_n\}\subset M_{3,2}$, then in particular \
$\{x_n\}\subset U_2\subset R_2$. Recall that $U_1\cup U_2=R_1\cup
R_2$. Therefore, \ $\dist(x_{n}, (U_1\cup U_2)^c)\ge \dist(x_{n},
R_2^c)$. Thus, $\lim_n \dist(x_{n}, R_2^c) =0$ and this implies
that $\lim_n g_2(x_n)=\gamma_2$.

\end{itemize}

\noindent Now, let us take any sequence $\{x_n\}\subset M_{3,i}$
such that $\lim_n g_i(x_n)=\gamma_i$. Consider every $x_n$ as an
element of $X^{**}$ and denote by $\bold{x_n}$ its restriction to
$[g_1,g_2,g_3]$. Recall that \ $D(x_n)\in S^*\cap [g_i,g_3]$ \ for
every $n\in \mathbb N$. Then, the sequence of restrictions
$\{\bold{x_n}\}\subset [g_1,g_2,g_3]^*$ satisfies that
\begin{align*}
\qquad 1&=|x_n|\ge|\bold{x_n}|\ge |\bold{x_n}|_{[g_i,g_3]}|=
\max\{\bold{x_n}(h):\ h\in S^*\cap[g_i,g_3]\}\\ &\ge
\bold{x_n}(D(x_n))=D(x_n)(x_n)=1,\end{align*}

\noindent for every $n\in \mathbb{N}$. Thus, there is a
subsequence \ $\{\bold{x_{n_j}}\}$ \ converging to an element \
$T\in [g_1,g_2,g_3]^*$ \ with \ $|T|=|T|_{[g_i,g_3]}|=1$. Since
$\lim_jg_i(x_{n_j})=\lim_{j}\bold{x_{n_j}}(g_i)=\gamma_i$\,, we
have that \ $T(g_i)=\gamma_i$ \ and this implies that \
$T|_{[g_i,g_3]}\in F_{3,i}'$ \ and \ $T\in F_{3,i}''$.
Furthermore, if $g\in N_{3,i}'$ and $T(g)=1$, then
$\lim_j\bold{x_{n_j}}(g)=1$. In addition,
$T(g_3)=\lim_j\bold{x_{n_j}}(g_3)=\lim_j g_3(x_{n_j})\ge \gamma_3$
because $\{x_{n_j}\}\subset U_3'$. Then, from condition
\eqref{3-3} if $i=1$ and condition \eqref{3-2} if $i=2$, we deduce
that $T(g_3)>\gamma_3$. Finally, let us check that
$T(g_{s})=\lim_j\bold{x_{n_j}}(g_{s})\le \gamma_{s}$\,, where
$s\in \{1,2\}$ and $s\not= i$\,:
\begin{itemize}
\item if $i=1$, the sequence $\{x_{n_j}\}\subset M_{3,1}$ and thus
$\{x_{n_j}\}\subset (U_1\cap U_3')\setminus U_2$. In particular \
$\{x_{n_j}\}\subset U_3'$ \ and  \ $g_1(x_{n_j})<
\gamma_{1,2}'<\gamma_{1,1}$ \ for every $j\in \mathbb N$.
Therefore, if \ $x_{n_j}\not\in U_2$ for all $j$, we must have
$\bold{x_{n_j}}(g_2)=g_2(x_{n_j})\le \gamma_2$ \ for every $j\in
\mathbb N$.
 \item if $i=2$, the sequence $\{x_{n_j}\}\subset
M_{3,2}$ and thus $\{x_{n_j}\}\subset (U_2\cap U_3')\setminus
U_1$. In particular \ $x_{n_j}\not\in U_1=R_1$\,, for every $j\in
\mathbb N$ \ and this implies \
$\bold{x_{n_j}}(g_1)=g_1(x_{n_j})\le \gamma_1$ \ for every $j\in
\mathbb N$.
\end{itemize}

\medskip

\item Consider a sequence $\{x_n\}\subset M_{3,1,2}$, such that
$\lim_n \dist(x_n,(U_1\cup U_2)^c)=0$. Then, it easily follows
that \ $\lim_n g_i(x_n)=\gamma_i$ \ for $i=1,2$. Indeed, $U_1\cup
U_2=R_1\cup R_2$ and then $\dist(x_n, (R_1\cup R_2)^c)\ge
\dist(x_n, R_i^c)$ for every $n\in \mathbb N$ and $i=1,2$. Hence
$\lim_n\dist(x_n, R_i^c)=0$. Since $\{x_n\}\subset R_i$,\  for
$i=1,2$, \ we obtain that $\lim_ng_i(x_n)=\gamma_i$, \ for
$i=1,2$.

\noindent Now, let us take any sequence $\{x_n\}\subset M_{3,1,2}$
such that $\lim_n g_i(x_n)=\gamma_i$, for every $i=1,2$. Consider
every $x_n$ as an element of $X^{**}$ and denote by $\bold{x_n}$
its restriction to $[g_1,g_2,g_3]$. Then, the sequence of
restrictions $\{\bold{x_n}\}\subset [g_1,g_2,g_3]^*$ satisfies
that
\begin{equation*}
\quad \qquad  1=|x_n|\ge|\bold{x_n}|=\max\{\bold{x_n}(h): \,h\in
S^*\cap[g_1,g_2,g_3]\}  \ge
\bold{x_n}(D(x_n))=D(x_n)(x_n)=1,\end{equation*}

\noindent for every $n\in \mathbb{N}$. Thus, there is a
subsequence \ $\{\bold{x_{n_j}}\}$ \ converging to an element \
$T\in [g_1,g_2,g_3]^*$ \ with \ $|T|=1$. Since \
$\lim_jg_i(x_{n_j})=\lim_{j}\bold{x_{n_j}}(g_i)=\gamma_i$ \ for \
$i=1,2$, \ then \ $T(g_i)=\gamma_i$ \ for \  $i=1,2$,\  and this
implies that $T\in F_{3,1,2}'$. Furthermore, if $g\in N_{3,1,2}'$
and $T(g)=1$, then $\lim_j\bold{x_{n_j}}(g)=1$. In addition,
$T(g_3)=\lim_j\bold{x_{n_j}}(g_3)=\lim_j g_3(x_{n_j})\ge \gamma_3$
because $\{x_{n_j}\}\subset U_3'$. Then, from condition
\eqref{3-1}, we deduce that $T(g_3)>\gamma_3$.
\end{enumerate}

\medskip

 \noindent Let us define, for $i=1,2$,
\begin{equation*}
F_{3,i}=\{T\in F_{3,i}'':\text{ there is } \ \{x_n\}\subset
M_{3,i} \, \text{ with }\,   \lim_n\bold{x_n}(g_i)=\gamma_i\,,
\text{ and }
 \lim_n\bold{x_n}=T\},
\end{equation*}
\begin{align*}
 F_{3,1,2}=\{T\in F_{3,1,2}':\text{ there is } \ \{x_n\}\subset
M_{3,1,2} \, \text{ with }\,   \lim_n\bold{x_n}(g_1)=\gamma_1\,, &
\ \lim_n\bold{x_n}(g_2)=\gamma_2  \\ & \quad \text{ and } \
 \lim_n\bold{x_n}=T\},
\end{align*}
and \begin{equation} F_3=F_{3,1}\cup F_{3,2}\cup F_{3,1,2}.
\end{equation}
Select a real number \ $\gamma_{3}'$ \  satisfying  \ $\gamma_3<
\gamma_{3}'<\min\{T(g_3):\ T\in F_3\}$ (recall that $F_3$ is
finite), and define,
\begin{equation*}
N_{3,i}=\{g\in N_{3,i}': \text{ there is } T\in  F_{3,i} \, \text{
with }\, T(g)=1 \}, \  \text{ for } i=1,2,
\end{equation*}
\begin{equation*}
 N_{3,1,2}=\{g\in  N_{3,1,2}':\text{ there is }  T \in
F_{3,1,2} \, \text{ with }\, T(g)=1\},
\end{equation*}
and $N_3=N_{3,1}\cup N_{3,2}\cup N_{3,1,2}$. Let us prove the
following Fact.

\medskip

\begin{fact} \label{3}
\begin{enumerate} \item There are numbers \ $0<t_3<l_3<1$ \ such that for every $g\in N_3$\,,
the slices $$O_g:=\{x\in S:\ g(x)>l_3\}  \  \text{ and } \
B_g:=\{x\in S:\ g(x)>t_3\}$$ satisfy that
\begin{align}\label{inclusion} O_g\subset B_g & \subset \{x\in S: \
g_1(x)<\gamma_{1,2}'\,,\ g_2(x)<\gamma_{2,1}'\,, \ g_3(x)>\gamma_{3}'\} \ \text{ and }\\
\label{intersection} &\dist(B_g,B_{g'})>0, \text{ whenever }
g,g'\in N_3\,, \ g\not=g'.
\end{align}
\item There are  numbers \  $\gamma_{1,2}\in
(\gamma_1,\gamma_{1,2}')$ \ and \  $\gamma_{2,1}\in
(\gamma_2,\gamma_{2,1}')$ \ such that if  $x\in M_3$\,, \
$g_1(x)<\gamma_{1,2}$ \ and  \ $g_2(x)<\gamma_{2,1}$\,, \ then \
$x\in O_g$\,, for some $g\in N_3$.
\end{enumerate}
\end{fact}

\medskip

\noindent{\bf Proof of Fact \ref{3}. } (1) First, if $X$ is
reflexive, we know that for every \ $g\in N_3$ \ there is \
$x_g\in S$ \ such that \ $D(x_g)=g$. Let us study the three
possible cases:
\begin{itemize}
\item If $g\in F_{3,1}$\,, denote by $\mathbf{x_g}$ the
restriction of $x_g$ to $[g_1,g_2,g_3]$. Since $\mathbf{x_g}(g)=1$
and $|\cdot|^*$ is G\^{a}teaux smooth, then $\bold{x_g}=T$ \ for
some $T\in F_{3,1}$. This implies that
$\bold{x_g}(g_1)=\gamma_1<\gamma_{1,2}'$\,, \
$\bold{x_g}(g_3)>\gamma_{3}'$ \ and \ $\bold{x_g}(g_2)\le
\gamma_2< \gamma_{2,1}'$.  Hence, $x_g\in \{x\in S: \
g_1(x)<\gamma_{1,2}'\,,\ g_2(x)>\gamma_{2,1}' \text{ and }
g_3(x)>\gamma_{3}'\}$.

\medskip

\item If $g\in F_{3,2}$\,, denote by $\mathbf{x_g}$ the
restriction of $x_g$ to $[g_1,g_2,g_3]$. Since $\mathbf{x_g}(g)=1$
and $|\cdot|^*$ is G\^{a}teaux smooth, then \ $\bold{x_g}=T$ \ for
some $T\in F_{3,2}$. This implies that
$\bold{x_g}(g_2)=\gamma_2<\gamma_{2,1}'$\,, \
$\bold{x_g}(g_3)>\gamma_{3}'$ \ and \ $\bold{x_g}(g_1)\le
\gamma_1< \gamma_{1,2}'$.  Hence, $x_g\in \{x\in S: \
g_1(x)<\gamma_{1,2}'\,,\ g_2(x)>\gamma_{2,1}' \text{ and }
g_3(x)>\gamma_{3}'\}$.

\medskip

\item  If $g\in F_{3,1,2}$\,,  denote by $\mathbf{x_g}$ the
restriction of $x_g$ to $[g_1,g_2,g_3]$. Since $\mathbf{x_g}(g)=1$
and $|\cdot|^*$ is G\^{a}teaux smooth, then $\bold{x_g}=T$ \ for
some $T\in F_{3,1,2}$. This implies that
$\bold{x_g}(g_1)=\gamma_1<\gamma_{1,2}'$\,, \
$\bold{x_g}(g_2)=\gamma_2<\gamma_{2,1}'$ \ and \ $\bold{x_g}(g_3)>
\gamma_{3}'$.  Hence, $x_g\in \{x\in S: \ g_1(x)<\gamma_{1,2}',\
g_2(x)>\gamma_{2,1}' \text{ and } g_3(x)>\gamma_{3}'\}$.
\end{itemize}

\medskip

 \noindent Now, since  the norm $|\cdot|$ is LUR
and $D(x_g)=g$, the functional $g$ strongly exposes $S$ at the
point $x_g$ for every $g\in N_3$. Since $N_3$ is finite, we can
hence obtain real numbers $0<t_3<l_3<1$ and slices $O_g$ and
$B_g$\,, \ for every $g\in N_3$\,, satisfying conditions
\eqref{inclusion} and \eqref{intersection}.

\medskip

\noindent Now consider a non reflexive Banach space $X$. Let us
first prove \eqref{inclusion}. Assume, on the contrary, that there
is a point $g \in N_3$ and there is a sequence $\{y_n\}\subset S$
satisfying $g(y_n)>1-\frac1n$ and such that $g_1(y_n)\ge
\gamma_{1,2}' $, or $g_2(y_n)\ge \gamma_{2,1}'$, or $g_3(y_n)\le
\gamma_3'$\,, for every $n\in \mathbb N$.
 If
$g\in N_{3}$  \ there is a sequence \ $\{x_n\}\subset M_{3}$ \
with \ $\lim_n g_i(x_n)\le \gamma_i$\,,\ for $i=1,2$, \ $\lim_n
g_3(x_n)>\gamma_3'$ \  and \ $\lim_n g(x_n)=1$. In particular,
\begin{equation*}
\frac{g(x_n)+1-\frac1n}{2}\le g\left(\frac{x_n+y_n}{2}\right)\le
\left|\frac{x_n+y_n}{2}\right|\le 1,
\end{equation*}
and thus $\lim_n\left|\frac{x_n+y_n}{2}\right|= 1$. Recall that in
the non reflexive case, the norm $|\cdot|$ is WUR, and  then
$x_n-y_n\xrightarrow{\omega} 0$\ (weaky converges to zero). This
last assertion gives a contradiction since we have
$\limsup_ng_1(x_n-y_n)\le \gamma_1-\gamma_{1,2}'<0$ \ or \
$\limsup_ng_2(x_n-y_n)\le \gamma_2-\gamma_{2,1}'<0$ \ or \
$\liminf_n g_3(x_n-y_n)\ge \lim_n g_3(x_n)-\gamma_3'>0$. Therefore
we can find real numbers $0<t_2<l_2<1$ and  slices $O_g$ and $B_g$
\ for every $g\in N_3$\,, satisfying condition \eqref{inclusion}.
The proof of \eqref{intersection} is the same as the one given in
Fact \ref{fact 1}, where the only property we need is the strict
convexity of $|\cdot|^*$.

\medskip

 \noindent (2)  Assume, on the contrary, that for every $n\in \mathbb N$, there is $x_n\in
M_3$ with $g_i(x_n)\le \gamma_i+\frac1n$, for $i=1,2$ \ and
$\{x_n:\ n \in \mathbb N\} \cap (\cup_{g\in N_3}O_g)=\emptyset$.
Then, there is a subsequence of $\{x_n\}$, which we keep denoting
by $\{x_n\}$, such that  $\{x_n\}\subset M_{3,1}$, \ or \
$\{x_n\}\subset M_{3,2}$, \ or \ $\{x_n\}\subset M_{3,1,2}$. In
the first case, $\lim_ng_1(x_n)=\gamma_1$. In the second case,
$\lim_ng_2(x_n)=\gamma_2$. In the third case,
$\lim_ng_i(x_n)=\gamma_i$\,, \ for every $i=1,2$. From the
definition of $F_3$ and $N_3$ and the comments preceding Fact
\ref{3}, we know that there is a subsequence $\{x_{n_j}\}$ and
$g\in N_3$ satisfying that $\lim_jg(x_{n_j})=1$, which is a
contradiction. This finishes the proof of Fact \ref{3}. $\Box$

\bigskip

If $R_3\cap (U_1\cup U_2)=\emptyset$ we may select as
$\gamma_{1,2}$ any number in $(\gamma_1,\gamma_{1,2}')$ and
$\gamma_{2,1}$ any number in $(\gamma_2,\gamma_{2,1}')$.

\medskip

Now we define $h_3$ as follows:
\begin{align*}
h_3&:\, S^+\longrightarrow \mathbb
R\\h_3&=\varphi_3(g_3)\,\phi_{2,1}(g_2)\,\phi_{1,2}(g_1),
\end{align*}
where $\varphi_3$,\  $\phi_{2,1}$ and $\phi_{1,2}$ are $C^\infty$
functions on $\mathbb R$ satisfying that
\begin{align*}
\varphi_3(t)&=0 \ \ \text{ if } t \le {\gamma_3}\\
\varphi_3(1)&=1\\
\varphi_3'(t)&>0 \ \ \text{ if } t>{\gamma_3}
\end{align*}
and
\begin{align*}
\phi_{1,2}(t) & =1 \ \ \text{ if } \ t\le
\textstyle{\frac{\gamma_1+\gamma_{1,2}}{2}} \qquad & \qquad
\phi_{2,1}(t) & =1 \ \ \text{ if } \ t\le
{\textstyle{\frac{\gamma_2+\gamma_{2,1}}{2}}}
\\
\phi_{1,2}(t) & =0 \ \ \text { if } \ t\ge \gamma_{1,2} \qquad &
\qquad \phi_{2,1}(t) & =0 \ \ \text{ if } \ t\ge \gamma_{2,1}
\\
\phi_{1,2}'(t) & <0 \ \ \text{ if } \ t \in
 \bigl(\textstyle{\frac{\gamma_1+\gamma_{1,2}}{2}}, \,
  \gamma_{1,2}\bigr) \qquad &  \qquad \phi_{2,1}'(t) & <0 \ \ \text{ if
}\ t\in \bigl( \textstyle{\frac{\gamma_2+\gamma_{2,1}}{2}}, \,
\gamma_{2,1} \bigr)
\end{align*}

\noindent Clearly the interior of the support of $h_3$ is the set
\begin{equation*}
U_3=\{x\in S^+:\ g_1(x)<\gamma_{1,2}\,, \ g_2(x)<\gamma_{2,1} \
\text{ and } \ g_3(x)>\gamma_3\}.\end{equation*} Select one point
$x_3\in U_3$, a real number $a_3\in \mathbb R^*$ with
$|a_3-F(x_3)|<\varepsilon$ \ and define the auxiliary function
\begin{align*}
&r_3:S^+\longrightarrow \mathbb R,\\
&r_3=s_3g_3+(1-s_3g_3(x_3)),\notag \end{align*} where we have
selected $s_3$ so that $s_3a_3>0$ and $|s_3|$ is small enough so
that the oscilation of $r_3$ on $U_3$ is less than \
$\frac{\varepsilon}{\,|a_3|}$. \, Notice that $r_3(x_3)=1$.

\medskip

\begin{figure}
  \includegraphics[width=14cm]{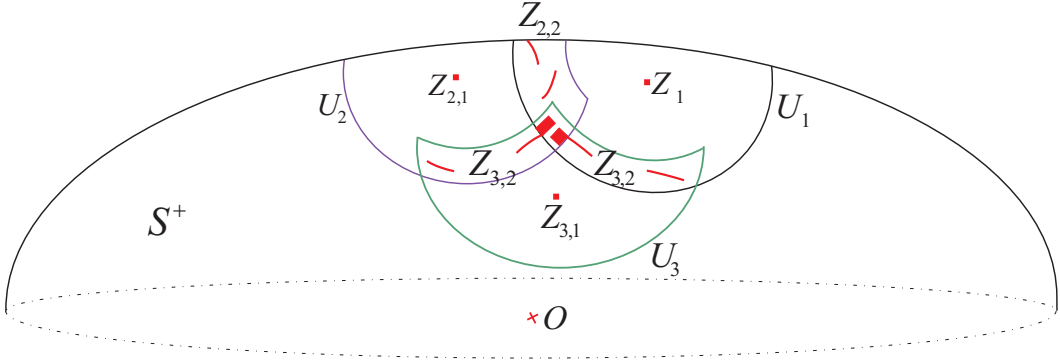}\\
  \caption{Case $n=3$: the decomposition of $Z_3$.}\label{dibujo n=3}
\end{figure}

\noindent Let us study the critical points $Z_3$ of the $C^p$
smooth function
\begin{align*} &{\bf H_3}: U_1\cup U_2\cup U_3\longrightarrow \mathbb R,\\ \notag
&{\bf H_3}=\frac{a_1r_1h_1+a_2r_2h_2+a_3r_3h_3}{h_1+h_2+h_3}.
\end{align*}
Let us prove that $Z_3:=\{x\in U_1\cup U_2\cup U_3:\,H'_3(x)=0
\text{ on } T_x\}$ can be included in a finite number of disjoint
slices within   $U_1\cup U_2\cup U_3$ by splitting it conveniently
into the (already defined) $Z_1$, $Z_2$ and up to four more
disjoint sets within $U_3$, as the Figure \ref{dibujo n=3}
suggests.

\medskip
\noindent The function ${\bf H_3'}$ can be written as ${\bf
H_3'}=\sigma_{3,1} \bold {g_1}+\sigma_{3,2} \bold
{g_2}+\sigma_{3,3} \bold {g_3}$,\ where \ $\sigma_{3,i}$ \ are
continuous and real functions on $U_1\cup U_2 \cup U_3$ \ and \
$\bold{g_i}$ \ denotes the restriction \ $g_i|_{T_x}$, \
$i=1,2,3$, whenever we evaluate $\bold{H_3'}(x)$ on $T_x$.

\medskip

\noindent Clearly $\bold{H_3}$ and $\bold{H_3'}$ restricted to
$(U_1\cup U_2)\setminus U_3$ coincide with $\bold{H_2}$ and
$\bold{H_2'}$ respectively. Then, $Z_3\setminus U_3=Z_3\setminus
\overline{U}_3=Z_2$. Let us study the set $Z_3\cap U_3$. \noindent
First, if $x\in U_3\setminus (U_1\cup U_2)$, \ then \
$H_3(x)=a_3r_3(x)$ \ and \ $H_3'(x)=a_3r_3'(x)=a_3s_3g_3$.
Therefore $\bold{H_3'}(x)=a_3s_3g_3|_{T_x}\equiv 0$ iff
$D(x)=g_3$. If the point $z_3\in U_3\setminus (U_1 \cup U_2)$
then, $H_3$ has exactly one critical point in $U_3\setminus
(U_1\cup U_2)$; in this case, since $g_1(z_3)\not=\gamma_1$\ and \
$g_2(z_3)\not={\gamma_2}$\,, the point $z_3$ actually belongs to
$U_3\setminus (\overline{U}_1 \cup \overline{U}_2)$.

\medskip

\noindent Now, let us study the critical points of ${\bf H_3}$ in
$U_3\cap(U_1\cup {U_2})$. If we define \
$\Lambda_2=\frac{h_1+h_2}{h_1+h_2+h_3}$, \ then we can rewrite \
$\bold{H_3}$ \ in \ $U_3\cap(U_1\cup U_2)$ \ as
\begin{equation*}
\bold{H_3}=\frac{a_1r_1h_1+a_2r_2h_2}{h_1+h_2}\,\cdot
\frac{h_1+h_2}{h_1+h_2+h_3}+\frac{a_3r_3h_3}{h_1+h_2+h_3}=\bold{H_2}\
\Lambda_2 +a_3r_3(1-\Lambda_2),
\end{equation*}

\noindent and \ $$\bold{{H_3}'}=\bold{{H_2}'}\Lambda_2
+a_3s_3(1-\Lambda_2)\bold{g_3}+(\bold{H_2}-a_3r_3)\Lambda'_2.$$ By
computing $\Lambda'_2$, we obtain
$\Lambda'_2=\xi_{2,1}\bold{g_1}+\xi_{2,2}\bold{g_2}+\xi_{2,3}
\bold{g_3}$, where the coefficients $\xi_{2,1},\ \xi_{2,2} $ \ and
\ $\xi_{2,3}$ are continuous functions of the following form:
 \begin{align}\label{derivada de Lambda}
\xi_{2,1}&=\frac{-\varphi_3(g_3)\phi_{2,1}(g_2)\phi_{1,2}'(g_1)(h_1+h_2)
+h_3\varphi_1'(g_1)+h_3\varphi_2(g_2)\phi_{1,1}'(g_1)}{(h_1+h_2+h_3)^2},\\
\xi_{2,2}&=
\frac{-\varphi_3(g_3)\phi_{2,1}'(g_2)\phi_{1,2}(g_1)(h_1+h_2)+h_3\varphi_2'(g_2)\phi_{1,1}(g_1)}{(h_1+h_2+h_3)^2},\,
\notag\\
\xi_{2,3}&=
\frac{-\varphi_3'(g_3)\,\phi_{2,1}(g_2)\,\phi_{1,2}(g_1)(h_1+h_2)}{(h_1+h_2+h_3)^2}.\,\notag
\end{align}

 \noindent
Since $g_1(x)< \gamma_{1,2}<\frac{\gamma_1+\gamma_{1,1}}{2}$ \ for
every $x\in U_3$\,, \  we have that $\phi_{1,1}'(g_1(x))=0$ \ for
every $x\in U_3$, and we can drop the term
$h_3\varphi_2(g_2)\phi_{1,1}'(g_1)$ in the above expression of
$\xi_{2,1}$. Thus, if $x\in U_3\cap(U_1\cup U_2)$, the
coefficients $\sigma_{3,1}\,, \  \sigma_{3,2}\,, \  \sigma_{3,3} $
\ for $\bold{H_3'}$ have the following form,
\begin{align*}
\sigma_{3,1}&=\sigma_{2,1}\Lambda_2+(\bold{H_2}-a_3r_3)\xi_{2,1}\\
\sigma_{3,2}&=\sigma_{2,2}\Lambda_2+(\bold{H_2}-a_3r_3)\xi_{2,2}\\
\sigma_{3,3}&=a_3s_3(1-\Lambda_2)+(\bold{H_2}-a_3r_3)\xi_{2,3},\\
\end{align*}
\noindent where  \ $a_3s_3>0$, \ $\Lambda_2>0$, \ $1-\Lambda_2>0$,
\ $\xi_{2,1}\ge 0$, \ $\xi_{2,2}\ge 0$,  \ $\xi_{2,1}+\xi_{2,2}>
0$ \ and \ $\xi_{2,3}<0$ \ on \ $U_3\cap (U_1\cup U_2)$.
Therefore, if $H_2-a_3r_3\le 0$, the coefficient $\sigma_{3,3}>0$.
When $H_2-a_3r_3 \ge 0$ and $\sigma_{2,2}>0$, we have that
$\sigma_{3,2}>0$. Finally, when $H_2-a_3r_3\ge 0$ and
$\sigma_{2,1}> 0$, we have $\sigma_{3,1}>0$ (recall that for every
$x\in U_1\cup U_2$, there is $j\in\{1,2\}$ such that
$\sigma_{2,j}(x)>0$). Since the vectors $\{g_1,g_2,g_3\}$ are
lineally independent we get that, if $\bold{H_3'}(x)= 0$ for some
$x\in U_3\cap (U_1 \cup U_2)$, then there necessarily exists
$\varrho\not=0$ such that $D(x)=\varrho(\sigma_{3,1}(x)g_1
+\sigma_{3,2}(x)g_2+ \sigma_{3,3}(x)g_3),$ \ that is
$D(x)\in[g_1,\, g_2,\,g_3].$ \

\medskip

\noindent In fact we can be more accurate and obtain that if $x\in
(U_3\cap U_2)\setminus U_1$ \ and \ ${\bf H_3'}(x)=0$ then
$D(x)\in [g_2,g_3]$. Indeed,  in  step 2 we proved that
$\sigma_{2,1}=0$ in $U_2\setminus U_1$. Moreover, the functions
$\varphi_1(g_1)$, $\phi_{1,1}(g_1)$ and $\phi_{1,2}(g_1)$ are
constant outside $U_1$, thus their derivatives vanish outside
$U_1$. This implies $\xi_{2,1}=0$ and consequently
$\sigma_{3,1}=0$ in $(U_3\cap U_2)\setminus U_1$. Similarly, if \
$x\in (U_3\cap U_1)\setminus U_2$ \ and ${\bf H_3'}(x)=0$,\ then
$D(x)\in [g_1,g_3]$. Indeed, from step 2 we know that
$\sigma_{2,2}=0$ on $U_1\setminus U_2$. Moreover, the function
$\varphi_2'(g_2)\phi_{1,1}(g_1)$ vanishes outside $U_2$. In
addition, if $x\in (U_3\cap U_1)\setminus U_2$ then
$g_1(x)<\gamma_{1,2}<\gamma_{1,1}$ and hence $g_2(x)\le \gamma_2$.
Thus $\phi_{2,1}'(g_2(x))=0$, which implies $\xi_{2,2}(x)=0$.
Consequently \ $\sigma_{3,2}(x)=0$  \ if \ $x\in(U_3\cap
U_1)\setminus U_2$.

\medskip

 \noindent Define the sets
\begin{align*}
Z_{3,1}&=\begin{cases}\{z_3\}, & \text{ if } \ z_3\in U_3\setminus
(\overline{U}_1\cup \overline{U}_2)\\ \emptyset, & \text{
otherwise }
\end{cases}\\
Z_{3,2}&=Z_3\cap U_3\cap (U_1\cup U_2).
\end{align*}
Now, let us check that $Z_{3,2}\subset\cup_{g\in N_3}O_g$. Indeed,
if $x\in Z_{3,2}$\,, then $x\in (U_1\cup U_2)\cap U_3$. Now,
\begin{itemize}
\item if $x\in (U_1\cap U_3)\setminus U_2$\,, then  $D(x)\in
[g_1,g_3]$. Since \ $ (U_1\cap U_3)\setminus U_2\subset (U_1\cap
U_3')\setminus U_2$ \ we can deduce that \ $x\in M_{3,1}\subset
M_3$. \item If \ $x\in (U_2\cap U_3)\setminus U_1$, then $D(x)\in
[g_2,g_3]$.  Since  $ (U_2\cap U_3)\setminus U_1\subset (U_2\cap
U_3')\setminus U_1$ we can deduce that $x\in M_{3,2}\subset M_3$.
\item If $x\in U_1\cap U_2 \cap U_3$, then $D(x)\in
[g_1,g_2,g_3]$. Since  $ U_1\cap U_2\cap U_3\subset U_1\cap
U_2\cap U_3'$ we can deduce that  $x\in M_{3,1,2}\subset M_3$.
\end{itemize}
Finally, since $x\in U_3$\,,\ we have that $g_1(x)<\gamma_{1,2}$
and $g_2(x)<\gamma_{2,1}$. We apply Fact \ref{3}(2) to conclude
that there is $g\in N_3$ such that $x\in O_g$.

\medskip

In the case when $Z_{3,1}=\{z_3\}\not\in \cup_{g\in
N_3}\overline{O}_g$, we select if necessary, a larger   $t_3$\,, \
with $t_3<l_3$\,,
 so that \ $z_3\not\in  \cup_{g\in
N_3}\overline{B}_g$. Since the norm is LUR and $D(z_3)=g_3$ we may
select numbers $0<t_3'<l_3'<1$ and open slices, which are
neighborhoods of $z_3$ defined by
\begin{equation*}
O_{g_3}:=\{x\in S: g_3(x)>l_3'\} \quad \text{ and } \quad
B_{g_3}:=\{x\in S: g_3(x)>t_3'\},
\end{equation*}
satisfying  \ $O_{g_3}\subset B_{g_3}\subset \{x\in S: \
g_1(x)<\gamma_{1,2}',\ g_2(x)<\gamma_{2,1}', \ g_3(x)>\gamma_3'\}
$ \ and \ $\dist(B_{g_3}, B_g)>0$ for every $g\in N_3$. In this
case, we define  \ $\Gamma_3=N_3\cup\{g_3\}$.

\medskip

Now, if $Z_{3,1}=\{z_3\}\in \cup_{g\in N_3}{\overline{O}_g}$,\ we
 select, if necessary, a smaller constant $l_3$\,, with $0<t_3<l_3<1$,  so that
 $Z_{3,1}=\{z_3\}\in \cup_{g\in N_3}{O_g}$\,. In this case, and also when $Z_{3,1}=\emptyset$,
 we define $\Gamma_3=N_3$.

 Notice that, in any of the cases mentioned above, Fact \ref{3} clearly
 holds for the (possibly) newly selected real numbers $t_3$ and $l_3$.

\medskip

Then, the distance between  any two sets  $B_{g}$, \  $B_{g'}$, \
$g,g'\in
 \Gamma_1\cup \Gamma_2 \cup \Gamma_3$, \
 $g\not=g'$, \  is strictly positive. \
Moreover \ $Z_{3,1}\cup Z_{3,2}\subset \cup_{g\in \Gamma_3} O_g
\subset  \cup_{g\in \Gamma_3} B_g\subset U_3'\subset R_3$.
Therefore, $Z_3=Z_1\cup Z_2 \cup Z_{3,1}\cup Z_{3,2}\subset
\cup_{g\in \Gamma_1 \cup \Gamma_2 \cup \Gamma_3 }O_g \subset
\cup_{g\in \Gamma_1 \cup \Gamma_2 \cup \Gamma_3 }B_g \subset
U_1\cup U_2\cup U_3=R_1\cup R_2 \cup R_3$. Finally, recall that \
$\dist(B_g, R_3^c)>0$, \ for every \ $g\in \Gamma_3$ \ and \
$\dist(B_g\,,\, (U_1\cup U_2\cup U_3)^c)>0$ \ for every \ $g\in
\Gamma_1\cup \Gamma_2 \cup \Gamma_3$.

\medskip

 It is worth mentioning that, by combining all the results obtained
in the step $n=3$, we deduce that ${\bf H_3'}=\sigma_{3,1} \bold
{g_1}+\sigma_{3,2} \bold {g_2}+\sigma_{3,3} \bold {g_3}$,\ where \
$\sigma_{3,i}$  \ are continuous functions on \ $U_1\cup U_2 \cup
U_3$\,, \ $\sigma_{3,i}(x)=0$ \ whenever \  $x\not \in U_i$\,, \
and for every \ $x\in U_1\cup U_2 \cup U_3$ \ there is at least
one coefficient \ $\sigma_{3,j}(x)>0$.

\bigskip
\bigskip

\noindent $\bold{\diamond}$ Assume that, in the steps {\bf $j=2,
..., k$}, with $k\geq 2$, we have selected points $z_j\in S^+$ and
constants $\gamma_j\in (0,1)$,  with \
$g_{1}(z_j)\not=\gamma_{1}$\,,\,...,\,$g_{j-1}(z_j)\not=\gamma_{j-1}$\,,
\ $\{g_1,...,g_k:=D(z_k)\}$ \ linearly independent functionals
such that
 \begin{align}
S_j= \{x\in S:\ f_j(x)>\delta_j\}\subset\{x\in S:\
g_j(x)>{\gamma_j}\}\subset \{x\in S: \ f_j(x)>\delta_j^4\}=P_j,
 \end{align}
for all $j=2,...,k$, and \begin{equation*} \{T\in
[g_{i_1},...,g_{i_s},g_j]^*: \ g_{i_1}(x)=\gamma_{i_1}\,,\, ...
\,,g_{i_s}(x)=\gamma_{i_s}\,,\ g_j(x)=\gamma_j \ \text{ and } \
|T|=1 \}=\emptyset,
\end{equation*}
for every  $1\le i_1<...<i_s\le j-1,\ \text{ and } \  1\le s\le
j-1$, $2\leq j\leq k$. Assume we have defined the functions \
$h_j=\varphi_j(g_j)\,\phi_{j-1,1}(g_{j-1})\,\cdots\phi_{1,j-1}(g_1),$
\ where \ $\varphi_j$\,,\ $\phi_{j-1,1}$\,,\,...,\,$\phi_{1,j-1}$
are $C^\infty$ functions on $\mathbb R$ satisfying
\begin{align*}
\varphi_j(t)&=0 \ \ \text{ if } t \le {\gamma_j}\\
\varphi_j(1)&=1\\
\varphi_j'(t)&>0 \ \ \text{ if } t>{\gamma_j}
\end{align*}
and
\begin{align*}
\phi_{1,j-1}(t) & =1 \ \ \text{ if } \ \textstyle{ t\le
\frac{\gamma_1+\gamma_{1,j-1}}{2}},&....., & \ \ \phi_{j-1,1}(t)
=1 \ \ \text{ if } \ \textstyle{t\le
{\frac{\gamma_{j-1}+\gamma_{j-1,1}}{2}}}
\\
\phi_{1,j-1}(t) & =0 \ \ \text { if } \ t\ge
\gamma_{1,j-1},&....., & \ \ \phi_{j-1,1}(t)  =0 \ \ \text{ if } \
t\ge
 \gamma_{j-1,1}
\\
\phi_{1,j-1}'(t) & <0 \ \ \text{ if } \ \textstyle{t \in
 \bigl(\frac{\gamma_1+\gamma_{1,j-1}}{2}}, \, \gamma_{1,j-1}\bigr),&....., &
 \  \ \phi_{j-1,1}'(t) <0 \ \ \text{ if
}\ t\in \bigl( \textstyle{\frac{\gamma_{j-1}+\gamma_{j-1,1}}{2}},
\, \gamma_{j-1,1} \bigr),
\end{align*}
where \ $\gamma_1<\gamma_{1,j-1}\,,......,
\gamma_{j-1}<\gamma_{j-1,1}$, and $2\leq j\leq k$.

\medskip

\noindent The interior of the support of $h_j$ is the set
\begin{equation*}
U_j=\{x\in S:\
 g_1(x)<\gamma_{1,j-1}\,,..., \, g_{j-1}(x)<\gamma_{j-1,1} \  \text{ and } \
 g_j(x)>\gamma_j\}. \end{equation*}
Assume we have also defined \ the $C^p$ smooth functions \ $r_j$ \
and \ $\bold{H_j}$:
\begin{align*}
&r_j:S^+\longrightarrow \mathbb R, \qquad \quad & {\bf H_j}&:
U_1\cup
U_2\cup...\cup U_j\longrightarrow \mathbb R,\\
&r_j=s_j g_j+(1-s_jg_j(x_j)),\notag \qquad \quad &{\bf
H_j}&=\frac{\sum_{i=1}^j a_ir_ih_i}{\sum_{i=1}^j h_i},
\end{align*}
for $2\leq j\leq k$, where $x_j\in U_j$\, the numbers $a_j\,,\,
s_j\in \mathbb R^*$ satisfy that $|a_j-F(x_j)|<\varepsilon$, \
$s_j a_j>0$, and the oscillation of $r_j$ on $U_j$ is less than \
$\frac{\varepsilon}{\,|a_j|}$.

\noindent  Assume that for $2\le j \le k$ the set of critical
points \ $Z_j$ \ of \ ${\bf H_j}$ \  is a union of the form
$Z_j=Z_{j-1}\cup Z_{j,1}\cup Z_{j,2}\,,$ where $Z_{j-1}$ is the
set of critical points of ${\bf H_{j-1}}$, the sets
$Z_{j-1},\,Z_{j,1},\,Z_{j,2}$ are pairwise disjoint, \
$Z_{j}\subset D^{-1}([g_1,...,g_j])$,  \ $Z_{j-1}\subset
(U_1\cup...\cup U_{j-1})\setminus \overline{U}_j$ \ and \
$Z_{j,1}\cup Z_{j,2}\subset U_j$. Furthermore, assume that there
is an open subset $U_j'$ such that \ $$U_j\subset U_j'\subset
R_j:=\{x\in S:\ g_j(x)>\gamma_j \}$$ \ and $\dist(B_g\,,U_j')>0$,
\ for every $g\in \Gamma_1\cup...\cup \Gamma_{j-1}$,\ there is a
finite subset $\Gamma_j\subset S^*$ and open slices of $S$,
\begin{equation*}
B_g:=\{x\in S: g(x)>t_j\} \quad \text{ and } \quad O_g:=\{x\in S:
\ g(x) > l_j\}, \quad 0<t_j<l_j<1,
\end{equation*}
satisfying $B_g\subset U_j'$ \ for every $g\in \Gamma_j$,  \
$\dist(B_g,B_{g'})>0$ \ whenever $g,g'\in \Gamma_j$, \ $g\not=g'$
\ and there is $\gamma_j'\in (\gamma_j,1)$ such that
$$Z_{j,1}\cup Z_{j,2}\subset \cup_{g\in \Gamma_j} O_g\subset
\cup_{g\in \Gamma_j} B_g\subset U_j'\cap \{x\in S:\
g_j(x)>\gamma_j'\}\subset U_j.$$
\medskip

\noindent  Finally, assume that for $2\le j\le k$, \  ${\bf
H_j'}=\sigma_{j,1}\,\bold{g_1}+\cdots+\sigma_{j,j}\,\bold{g_j}$ \
on $U_1\cup...\cup U_j$\,, \ where \ $\sigma_{j,i}$ \ are
continuous functions on \ $U_1\cup...\cup U_j$, for every \
$i=1,....,j$.\ and assume that for every
 $x\in U_1 \cup \dots \cup U_j$ there is at least one
strictly positive coefficient \ $\sigma_{j,m}(x)$, \  and that if
\ $x\in (U_1\cup...\cup U_j)\setminus U_m$ , \ with \ $m\in
\{1,...,j\}$, \ then $\sigma_{j,m}(x)=0$.

\bigskip
\bigskip

\noindent $\bold{\diamond}$ Now, let us denote by $y_{k+1}\in S$
the point satisfying $f_{k+1}(y_{k+1})=1$.  If either
$\{g_1,...,g_k,f_{k+1}\}$ are lineally dependent or
$g_i(y_{k+1})=\gamma_i$ for some $i\in \{1,...,k\}$, we can use
the density of the norm attaining functionals (Bishop-Phelps
Theorem) and the continuity of $D$ to slightly modify $y_{k+1}$
and find $z_{k+1}\in S$ so that: $ g_i(z_{k+1})\not=\gamma_i$\,, \
for every $i=1,...,k$, \ $\{g_1\,,...,g_k\,,g_{k+1}:=D(z_{k+1})\}$
are l.i. and
 \begin{align*}
 \{x\in S:\ f_{k+1}(x)>\delta_{k+1}^2\}\subset\{x\in S:\
g_{k+1}(x)>{\nu_{k+1}}\}\subset \{x\in S: \
f_{k+1}(x)>\delta_{k+1}^3\},
 \end{align*}
for some $\nu_{k+1}\in (0,1)$. If $g_i(y_{k+1})\not=\gamma_i$ for
every $i\in \{1,...,k\}$ and $\{g_1,...,g_k,f_{k+1}\}$ are
lineally independent, we define $z_{k+1}=y_{k+1}$ and
$g_{k+1}=f_{k+1}$. Then we apply Lemma \ref{case N} to the l.i.
vectors, $\{g_1,...,g_{k+1}\}$ and the real numbers
$\gamma_1,....,\gamma_{k}$ and obtain $\gamma_{k+1}\in (0,1)$
close enough to $\nu_{k+1}$ \ so that
\begin{align*}
 S_{k+1}&=\{x\in S:\ f_{k+1}(x)>\delta_{k+1}\}\\ &\subset\{x\in S:\
g_{k+1}(x)>{\gamma_{k+1}}\}\subset \{x\in S: \
f_{k+1}(x)>\delta_{k+1}^4\}=P_{k+1}
\end{align*}
and
 \begin{equation}\label{empty k+1}
 \{T\in [g_{i_1}\,,...,\,g_{i_s}\,,g_{k+1}]^* :\,T(g_{i_1})=\gamma_{i_1}\,,...,\,T(g_{i_s})=\gamma_{i_s}\,,
T(g_{k+1})=\gamma_{k+1} \ \text{ and } \  |T|=1 \}=\emptyset
\end{equation} for every $1\le i_1<...<i_s\le k$ \ and  \ $ 1\le s\le
k$.

\medskip

\noindent Define $$R_{k+1}=\{x\in S: \
g_{k+1}(x)>\gamma_{k+1}\}.$$ Recall that $\cup_{g\in \Gamma_k}
B_g\subset U_k'\cap \{x\in S: \ g_k(x)>\gamma_k'\}$ \ and select \
$\gamma_{k,1}'\in (\gamma_k\,,\gamma_k')$. In addition, we select
numbers
\begin{equation}\label{half} \gamma_{k-1,2}'\in
\big(\gamma_{k-1}\,,\,\frac{\gamma_{k-1}+\gamma_{k-1,1}}{2}\big)\,,....,\,\gamma_{1,k}'\in
\big(\gamma_1\,,\,\frac{\gamma_1+\gamma_{1,k-1}}{2}\big),
\end{equation}
and define the open set
\begin{equation}U_{k+1}'=\{x\in S:\
 g_1(x)<\gamma_{1,k}'\,,..., \, g_{k}(x)<\gamma_{k,1}' \  \text{ and }\
 g_{k+1}(x)>\gamma_{k+1}\}.
\end{equation}
Notice that \ $\operatorname{dist}(B_g,\,U_{k+1}')>0$ \ for every
\ $g\in \Gamma_1\cup...\cup \Gamma_k $.

\medskip

 Assume that
$R_{k+1}\cap(U_1\cup...\cup U_k)\not=\emptyset$ and define, for
every $1\le i_1<...<i_s \le k$ \ and \ $1\le s\le k$\,, the set
\begin{align*}
M_{k+1,i_1,...,i_s}=\bigl\{x\in U_{i_1}\cap...\cap U_{i_s}\cap
U_{k+1}':\ x\not\in U_j \
\text{ for every } & \ j\in \{1,...,k\} \setminus \{i_1,...,i_s\}, \\
 & \text{ and  }  D(x)\in[g_{i_1},...,g_{i_s},g_{k+1}] \bigr\}.
\end{align*}
and $$M_{k+1}=\bigcup\{M_{k+1,i_1,...,i_s}: \ 1\le i_1<...<i_s \le
k \ \text{ and } \ 1\le s\le k\}.$$ In the case when
$M_{k+1}=\emptyset$ we select as $\gamma_{1,k}$ any point in
$(\gamma_1,\gamma_{1,k}')$,...., and $\gamma_{k,1}$ any point in
$(\gamma_k,\gamma_{k,1}')$.

Notice  that $U_1\cup ...\cup U_k=R_1\cup...\cup R_k$.  In the
case when \ $M_{k+1}\not=\emptyset$ \ and \ $\dist(M_{k+1}\,,
(U_1\cup...\cup U_k)^c)=\dist(M_{k+1}\,, (R_1\cup...\cup
R_k)^c)>0$, we can immediately find $\gamma_{1,k}\in
(\gamma_1\,,\gamma_{1,k}')$\,,....,$\gamma_{k\,,1}\in(\gamma_k\,,\gamma_{k,1}')$
with $M_{k+1}\subset \{x\in S: g_1(x)>\gamma_{1,k}\}\cup...\cup
\{x\in S: g_k(x)>\gamma_{k,1}\}$.

\medskip

In the case when $\dist(M_{k+1}\,, (U_1\cup...\cup U_k)^c)=0$ and
in order to find suitable positive numbers
$\gamma_{1,k}\,,...,\gamma_{k,1}$, we need to study the limits of
the sequences $\{x_n\} \subset M_{k+1}$ such that $\lim_n
\dist(x_n, (U_1\cup...\cup U_k)^c)=0$. Define, for every $1\le
i_1<...<i_s \le k$ and  $1\le s\le k$, the sets
\begin{align*}
F_{k+1,i_1,...,i_s}'&=\{T\in [g_{i_1},...,g_{i_s},g_{k+1}]^*: \
T(g_i)=\gamma_i \ \text{ for every } \ i\in\{i_1,...,i_s\} \text{
and } \ |T|=1\},\\ N_{k+1,i_1,...,i_s}'&=\{g\in S^*\cap
[g_{i_1},...,g_{i_s},g_{k+1}]:\ T(g)=1 \ \text{ for some } \ T\in
F_{k+1,i_1,...,i_s}'\}.
\end{align*}
Since the norm $|\cdot|^*$ is G\^{a}teaux smooth,  we can apply
Lemma \ref{vectoradicional} to the finite dimensional space
$[g_{i_1},...,g_{i_s},g_{k+1}]$ with the norm $|\cdot|^*$
restricted to this finite dimensional space, and deduce that the
cardinal of any of the sets $F_{k+1,i_1,...,i_s}'$ is at most two.
Moreover, since the norm is strictly convex, the cardinal of each
set $N_{k+1,i_1,...,i_s}'$ is at most two. Let us consider the
{\em norm-one} extensions to \ $[g_1\,,...,g_k\,,g_{k+1}]$ \ of
the elements of \ $F_{k+1,i_1,...,i_s}'$\,, that is,
\begin{equation*}
F_{k+1,i_1,...,i_s}''=\{T\in [g_1,...,g_{k+1}]^*:\
T|_{[g_{i_1},...,g_{i_s},g_{k+1}]}\in F_{k+1,i_1,...,i_s}' \
\text{ and } \ |T|=1 \}.
\end{equation*}
Since the norm $|\cdot|^*$ is G\^{a}teaux smooth, for every $G\in
F_{k+1,i_1,...,i_s}'$ there is a {\em unique norm-one } extension
$T$ defined on $[g_1\,,...,g_{k+1}]$. Thus, the cardinal of the
every set $F_{k+1,i_1,...,i_s}''$ is at most two. Therefore the
sets
\begin{align*} F_{k+1}'&=\bigcup \{F_{k+1,i_1,...,i_s}'': \ 1\le
i_1<...<i_s \le k \ \text{ and } \ 1\le s\le k\}
\\N_{k+1}'&=\bigcup\{N_{k+1,i_1,...,i_s}':\ 1\le i_1<...<i_s \le k
\ \text{ and } \ 1\le s\le k \}
\end{align*} are finite.  As a
consequence of equality \eqref{empty k+1}, we deduce that
$T(g_{k+1})\not=\gamma_{k+1}$ for every $T\in F_{k+1}'$. We can
restrict our study to the following kind of sequences: Fix \ $
1\le s\le k$ \ and \ $1\le i_1<...<i_s \le k$ \ and consider a
sequence $\{x_n\}\subset M_{k+1,i_1,...,i_s}$ \ such that
$\lim_n\dist(x_n, (U_1\cup ...\cup U_k)^c)=0$. Let us prove that
for every \ $i\in\{i_1,...,i_s\}$, \ $\lim_n g_i(x_n)=
\gamma_i$\,. Indeed, if $\{x_n\}\subset M_{k+1,i_1,...,i_s}$\,,
then in particular $\{x_n\}\subset U_i\subset R_i$ for every $i\in
\{i_1\,,...,i_s\}$. Recall that $U_1\cup...\cup U_k=R_1\cup
...\cup R_k$. Therefore, $\dist(x_n\,,(U_1\cup...\cup U_k)^c)\ge
\dist(x_n\,,R_i^c)$ \ for every $i\in \{i_1,...,i_s\}$. Thus,
$\lim_n\dist(x_n\,,R_i^c)=0$ \ for every $i\in \{i_1,...,i_s\}$.
Since $\{x_n\}\subset R_i$\,, this implies that
$\lim_ng_i(x_n)=\gamma_i$, \ for every $i\in \{i_1,...,i_s\}$.

\medskip

Now, let us take any sequence $\{x_n\} \subset
M_{k+1,i_1,...,i_s}$ such that $\lim_n g_i(x_n)=\gamma_i$\,, \ for
$i\in\{i_1,...,i_s\}$.  Consider every $x_n$ as an element of
$X^{**}$ and denote by $\bold{x_n}$ its restriction to
$[g_1,...,g_{k+1}]$. Recall that $D(x_n)\in S^*\cap
[g_{i_1},...,g_{i_s},g_{k+1}]$ for every $n\in \mathbb N$. Then,
the sequence of restrictions $\{\bold{x_n}\} \subset
[g_1,...,g_{k+1}]^*$ satisfies that
\begin{align*}
1=&|x_n| \ge |\bold{{x_n}}|\ge
|\bold{{x_n}}|_{[g_{i_1},...,g_{i_s},g_{k+1}]}|= \max
\{\bold{x_n}(h):\ h\in S^*\cap [g_{i_1},...,g_{i_s},g_{k+1}] \}
\\&\ge \bold{x_n}(D(x_n))=D(x_n)(x_n)=1,
\end{align*}
for every $n\in \mathbb N$. Thus, there is a subsequence
$\{\bold{x_{n_j}}\}$ converging to an element $T\in
[g_1,...,g_{k+1}]^*$ \ with \ $|T|=1$ \ and \
$|T|_{[g_{i_1},...,g_{i_s},g_{k+1}]}|=1$. Since \ $\lim_j
g_i(x_{n_j})= \gamma_i$ \ for every \ $i\in \{i_1,...,i_s\}$, \ we
have that $T(g_i)=\gamma_i$ \ for every \ $i\in \{i_1,...,i_s\}$.
\ This implies that $T|_{[g_{i_1},...,g_{i_s},g_{k+1}]}\in
F_{k+1,i_1,...,i_s}'$ \ and \ $T\in F_{k+1,i_1,...,i_s}''$.
Furthermore, \ if \ $g\in N_{k+1,i_1,...,i_s}'$ \ and \ $T(g)=1$,
\ then \ $\lim_j \bold{x_{n_j}}(g)=1$.  In addition,
$T(g_{k+1})=\lim_j
\bold{x_{n_j}}(g_{k+1})=\lim_jg_{k+1}(x_{n_j})\ge \gamma_{k+1}$
because $\{x_{n_j}\}\subset U_{k+1}'$. Then, from condition
\eqref{empty k+1}, we deduce that \ $T(g_{k+1})>\gamma_{k+1}$.
 Finally, let us check that $g_i(x_{n})\le
\gamma_i$ for every \ $i\in\{1,...,k\}\setminus \{i_1,...,i_s\}$ \
and \ $n\in \mathbb N$. Indeed, since $\{x_n\}\subset U_{k+1}'$,
twe have
$g_1(x)<\gamma_{1,k}'<\gamma_{1,i-1},...,g_{i-1}(x)<\gamma_{i-1,k+2-i}'<\gamma_{i-1,1}$.
Now, from the definition of $U_i$ and the fact that  $\{x_n:\ n\in
\mathbb N\}\cap U_i=\emptyset$, we deduce that $g_i(x_n)\le
\gamma_{i}$, for every $n\in \mathbb N$. Finally, if \ $T=\lim_j
\bold{x_{n_j}}$ \ in \ $[g_1,...,g_{k+1}]$, then \ $T(g_i)=\lim_j
\bold{x_{n_j}}(g_i)\le \gamma_i$,\  for every \
$i\in\{1,...,k\}\setminus \{i_1,...,i_s\}$.

\bigskip

\noindent Let us define, for every $1\le s\le k$ and $1\le
i_1<...<i_s\le k$, the sets
\begin{align*}
F_{k+1,i_1,...,i_s}=\{T\in F_{k+1,i_1,...,i_s}'':\  \text{ there
is }  \{x_n\} &\subset M_{k+1,i_1,...,i_s} ,  \text{ with }\
\lim_n \bold{x_n}(g_i)=\gamma_i, \\ & \quad  \text{ for } \
i\in\{i_1,...,i_s\} \ \text{ and } \ \lim_n\bold{x_n}=T\},
\end{align*}
\begin{equation*}
N_{k+1,i_1,...,i_s}=\{g\in N_{k+1,i_1,...,i_s}':\ \text{ there is
} T\in F_{k+1,i_1,...,i_s} \text{ with } T(g)=1\},
\end{equation*}
and
 \begin{align*}
F_{k+1}&=\bigcup \{F_{k+1,i_1,...,i_s}: \ 1\le s\le k \ \text{ and
} \  1\le i_1<...<i_s\le k\},\\
N_{k+1}&=\bigcup \{N_{k+1,i_1,...,i_s}: \ 1\le s\le k \ \text{ and
} \  1\le i_1<...<i_s\le k\},
 \end{align*}
which are all finite. Select a real number $\gamma_{k+1}'$
satisfying $\gamma_{k+1}<\gamma_{k+1}'<\min\{T(g_{k+1}): \ T\in
F_{k+1}\}$.

\medskip

\begin{fact}\label{k+1}
\begin{enumerate} \item There are numbers \ $0<t_{k+1}<l_{k+1}<1$ \ such that for every $g\in N_{k+1}$,
the slices $$O_g:=\{x\in S:\ g(x)>l_{k+1}\} \quad  \text{ and }
\quad B_g:=\{x\in S:\ g(x)>t_{k+1}\}$$ satisfy that
\begin{align}\label{inclusionk+1} O_g&\subset B_g\subset \{x\in S: \
g_1(x)<\gamma_{1,k}'\,,...,\ g_k(x)<\gamma_{k,1}'\,, \ g_{k+1}(x)>\gamma_{k+1}'\} \quad \text{ and }\\
\label{intersectionk+1} &\dist(B_g,B_{g'})>0, \text{ whenever }
g,g'\in N_{k+1}, \ g\not=g'.
\end{align}
\item There are  numbers\  $\gamma_{1,k}\in
(\gamma_1,\gamma_{1,k}')$,....,$\gamma_{k,1}\in
(\gamma_k,\gamma_{k,1}')$ such that if $x\in M_{k+1}$\,, \
$g_1(x)<\gamma_{1,k}$\,,....,\,$g_k(x)<\gamma_{k,1}$\,,\ then
$x\in O_g$\,, for some $g\in N_{k+1}$.
\end{enumerate}
\end{fact}

\medskip

\noindent{\bf Proof of Fact \ref{k+1}. } (1) First, if $X$ is
reflexive, we know that for every \ $g\in N_{k+1}$ \ there is \
$x_g\in S$ \ such that \ $D(x_g)=g$. There is $1\le s\le k$ and
$1\le i_1<...<i_s\le k$ such that $g\in F_{k+1,i_1,...,i_s}$.
Denote by $\mathbf{x_g}$ the restriction of $x_g$ to
$[g_1,...,g_{k+1}]$. Since $\mathbf{x_g}(g)=1$ and $|\cdot|^*$ is
G\^{a}teaux smooth, we have that $\bold{x_g}=T$ for some \ $T\in
F_{k+1,i_1,...,i_s}$. \ This implies that
$\bold{x_g}(g_i)=\gamma_i<\gamma_{i,k+1-i}'$\,, \ whenever
$i\in\{i_1,...,i_s\}$, \ $\bold{x_g}(g_{k+1})>\gamma_{k+1}'$ \ and
\ $\bold{x_g}(g_i)\le \gamma_i< \gamma_{i,k+1-i}'$\,, \ whenever
$i\in\{1,...,k\}\setminus\{i_1,...,i_s\}$. Hence, $x_g\in \{x\in
S: \ g_1(x)<\gamma_{1,k}'\,,...,\ g_k(x)<\gamma_{k,1}' \ \text{
and } \ g_{k+1}(x)>\gamma_{k+1}'\}$.

 \noindent Now, since  the norm $|\cdot|$ is LUR
and $D(x_g)=g$, the functional $g$ strongly exposes $S$ at the
point $x_g$ for every $g\in N_{k+1}$. Since $N_{k+1}$ is finite,
we can obtain real numbers $0<t_{k+1}<l_{k+1}<1$ and slices $O_g$
and $B_g$, \ for every $g\in N_{k+1}$, satisfying conditions
\eqref{inclusionk+1} and \eqref{intersectionk+1}.

\medskip

\noindent Now consider a non reflexive Banach space $X$. Let us
first prove \eqref{inclusionk+1}. Assume, on the contrary, that
there is a point \ $g \in N_{k+1}$ \ and there is a sequence \
$\{y_n\}\subset S$ \ satisfying \ $g(y_n)>1-\frac1n$ \ with  \
$g_1(y_n)\ge \gamma_{1,k}' $\,,...., or \ $g_k(y_n)\ge
\gamma_{k,1}'$, \ or \ $g_{k+1}(y_n)\le \gamma_{k+1}'$\,, \ for
every $n\in \mathbb N$.
 If
$g\in N_{k+1}$  \ there is a sequence \ $\{x_n\}\subset M_{k+1}$ \
with \ $\lim_n g_i(x_n)\le \gamma_i$\,,\ for every \
$i\in\{1,...,k\}$, \ $\lim_n g_{k+1}(x_n)>\gamma_{k+1}'$ \  and \
$\lim_n g(x_n)=1$. In particular,
\begin{equation*}
\frac{g(x_n)+1-\frac1n}{2}\le g\left(\frac{x_n+y_n}{2}\right)\le
\left|\frac{x_n+y_n}{2}\right|\le 1,
\end{equation*}
and thus $\lim_n\left|\frac{x_n+y_n}{2}\right|= 1$. Since in this
case the norm  $|\cdot|$  is WUR, we have that
$x_n-y_n\xrightarrow{\omega} 0$\ (weakly converges to zero). This
last assertion gives a contradiction since either \linebreak
$\limsup_ng_i(x_n-y_n)\le \gamma_i-\gamma_{i,k+1-i}'<0$ \ for some
$i\in\{1,...,k\}$ \ or \ $\liminf_n g_{k+1}(x_n-y_n)\ge \lim_n
g_{k+1}(x_n)-\gamma_{k+1}'>0$. Therefore, we can find real numbers
$0<t_{k+1}<l_{k+1}<1$ and  slices $O_g$ and $B_g$ \ for every
$g\in N_{k+1}$, satisfying condition \eqref{inclusionk+1}. The
proof of \eqref{intersectionk+1} is the same as the one given in
Fact \ref{k+1}, where the only property we need is the strict
convexity of $|\cdot|^*$.

\medskip

 \noindent (2)  Assume, on the contrary,
that for every $n\in \mathbb N$, there is $x_n\in M_{k+1}$ with
$g_i(x_n)\le \gamma_i+\frac1n$, for every $i\in\{1,...,k\}$ \ and
$\{x_n:\ n \in \mathbb N\} \cap (\cup_{g\in
N_{k+1}}O_g)=\emptyset$. Then there is a subsequence of $\{x_n\}$,
which we denote by $\{x_n\}$ as well, and there are numbers \
$1\le s\le k$ \ and \ $1\le i_1<...<i_s<k$ \ such that
$\{x_n\}\subset M_{k+1,i_1,...,i_s}$. In particular,  \
$\{x_n\}\subset U_i\subset R_i$  \ and then \ $g_i(x_n)>\gamma_i$
\ for every \ $i\in\{i_1\,,...,i_s\}$ \ and \ $n\in \mathbb N$.
Hence, \
 $\lim_ng_i(x_n)=\gamma_i$ \ for \ every \ $i\in\{i_1,...,i_s\}$. Since
$\{x_n\}\subset M_{k+1,i_1,...,i_s}$\,, from the comments
preceding Fact \ref{k+1}, we know that there is a subsequence
$\{x_{n_j}\}$ and $g\in N_{k+1,i_1,...,i_s}$ satisfying that
$\lim_jg(x_{n_j})=1$, which is a contradiction. This finishes the
proof of Fact \ref{k+1}. $\Box$

\medskip

If $R_{k+1}\cap (U_1\cup...\cup U_k)=\emptyset$ \ we may select as
$\gamma_{1,k}$ any number in $(\gamma_1, \gamma_{1,k}')$\,,....,
and \ $\gamma_{k,1}$ \ any number in \ $(\gamma_k,
\gamma_{k,1}')$.

\medskip

Now we define $h_{k+1}$,
\begin{align*}
h_{k+1} & : S^+\longrightarrow \mathbb R
\\
h_{k+1} &
=\varphi_{k+1}(g_{k+1})\,\phi_{k,1}(g_{k})\,\cdots\phi_{1,k}(g_1),\end{align*}
 with \ $\varphi_{k+1}$,\ $\phi_{k,1}$,\,...,\,$\phi_{1,k}$
$C^\infty$ functions on $\mathbb R$ satisfying
\begin{align*}
\varphi_{k+1}(t)&=0 \ \ \text{ if } t \le {\gamma_{k+1}}\\
\varphi_{k+1}(1)&=1\\
\varphi_{k+1}'(t)&>0 \ \ \text{ if } t>{\gamma_{k+1}}
\end{align*}
and
\begin{align*}
\phi_{1,k}(t) & =1 \ \ \text{ if } \ \textstyle{ t\le
\frac{\gamma_1+\gamma_{1,k}}{2}},&....., & \ \ \phi_{k,1}(t) =1 \
\ \text{ if } \ \textstyle{t\le
{\frac{\gamma_{k}+\gamma_{k,1}}{2}}}
\\
\phi_{1,k}(t) & =0 \ \ \text { if } \ t\ge \gamma_{1,k},&....., &
\ \ \phi_{k,1}(t)  =0 \ \ \text{ if } \ t\ge
 \gamma_{k,1}
\\
\phi_{1,k}'(t) & <0 \ \ \text{ if } \ \textstyle{t \in
 \bigl(\frac{\gamma_1+\gamma_{1,k}}{2}}, \, \gamma_{1,k}\bigr),&....., &
 \  \ \phi_{k,1}'(t) <0 \ \ \text{ if
}\ t\in \bigl( \textstyle{\frac{\gamma_{k}+\gamma_{k,1}}{2}}, \,
\gamma_{k,1} \bigr),
\end{align*}
\noindent Clearly the interior of the support of $h_{k+1}$ is the
set \begin{equation*} U_{k+1}=\{x\in S: \
g_1(x)<\gamma_{1,k}\,,...,\,g_{k}(x)<\gamma_{k,1} \ \text{ and } \
g_{k+1}(x)>\gamma_{k+1}\}.\end{equation*} Select one point
$x_{k+1}\in U_{k+1}$, a real number $a_{k+1}\in \mathbb R^*$ with
$|a_{k+1}-F(x_{k+1})|<\varepsilon$ and define the auxiliary
function
\begin{align*}
&r_{k+1}:S^+\longrightarrow \mathbb R,\\
&r_{k+1}=s_{k+1}g_{k+1}+(1-s_{k+1}g_{k+1}(x_{k+1})),\notag
\end{align*} where we have selected $s_{k+1}$ so that $s_{k+1}a_{k+1}>0$ and
$|s_{k+1}|$ is small enough so that the oscilation of $r_{k+1}$ on
$U_{k+1}$ is less than \ $\frac{\varepsilon}{\,|a_{k+1}|}$. \,
Notice that $r_{k+1}(x_{k+1})=1$.

\noindent Let us study the set of critical points $Z_{k+1}$ of the
$C^p$ smooth function
\begin{align*} &{\bf H_{k+1}}: U_1\cup ...\cup U_{k+1}\longrightarrow \mathbb R,\\ \notag
&{\bf H_{k+1}}=\frac{\sum_{i=1}^{k+1} a_ir_ih_i}{\sum_{i=1}^{k+1}
h_i}.
\end{align*}
Let us prove that $Z_{k+1}:=\{x\in U_1\cup ... \cup U_{k+1}: \,
H'_{k+1}(x)=0 \text{ on } T_x \}$  can be included in a finite
union of disjoint slices within  $U_1\cup ... \cup U_{k+1}$ by
splitting $Z_{k+1}$ conveniently into the (already defined)  set
$Z_k$ and a finite number of disjoint sets within $U_{k+1}$.

\medskip
\noindent It is straightforward to verify that ${\bf
H_{k+1}'}=\sigma_{k+1,1} \bold {g_1}+...+\sigma_{k+1,k+1}
\bold{g_{k+1}}$,\ where \ $\sigma_{k+1,i}$ \ are continuous
functions on $U_1\cup ... \cup U_{k+1}$ \ and \ $\bold{g_i}$ \
denotes the restriction \ $g_i|_{T_x}$, \ $i=1,...,k+1$, whenever
we evaluate ${\bf H_{k+1}'}(x)$.

\medskip

\noindent Clearly the restrictions of $\bold{H_{k+1}}$ and
$\bold{H_{k+1}'}$ to $(U_1\cup ...\cup U_k)\setminus U_{k+1}$
coincide with $\bold{H_k}$ and $\bold{H_k'}$ respectively. Then,
$Z_{k+1}\setminus U_{k+1}=Z_k=Z_{k+1}\setminus
\overline{U}_{k+1}$. Let us study the set $Z_{k+1}\cap U_{k+1}$.
\noindent First, if $x\in U_{k+1}\setminus (U_1\cup ...\cup
U_k)$,\ then \ $H_{k+1}(x)=a_{k+1}r_{k+1}(x)$ and
$H_{k+1}'(x)=a_{k+1}r_{k+1}'(x)$. Therefore
$\bold{H_{k+1}'}(x)=a_{k+1}s_{k+1}g_{k+1}|_{T_x}= 0$ \ iff \
$D(x)=g_{k+1}$. If the point $z_{k+1}\in U_{k+1}\setminus (U_1
\cup ...\cup U_k)$, then ${\bf H_{k+1}}$ has exactly one critical
point in $U_{k+1}\setminus (U_1\cup ...\cup U_k)$; in this case,
since $g_i(z_{k+1})\not=\gamma_i$ \ for every  \ $i=1,..,k$, \ the
point $z_{k+1}$ actually belongs to $U_{k+1}\setminus
(\overline{U}_1 \cup ... \cup \overline{U}_k)$.

\medskip

\noindent Now, let us study the critical points of ${\bf H_{k+1}}$
in $U_{k+1}\cap(U_1\cup ...\cup U_k)$. If we define \
$\Lambda_{k}=\frac{\sum_{i=1}^k h_i}{\sum_{i=1}^{k+1} h_i}$, \
then we can rewrite \ $\bold{H_{k+1}}$ \ on \ $U_{k+1}\cap(U_1\cup
...\cup U_k)$ \ as
\begin{equation*}
\bold{H_{k+1}}=\frac{\sum_{i=1}^{k}a_ir_ih_i}{\sum_{i=1}^k
h_i}\,\cdot \frac{\sum_{i=1}^k h_i}{\sum_{i=1}^{k+1}
h_i}+\frac{a_{k+1}r_{k+1}h_{k+1}}{\sum_{i=1}^{k+1}
h_i}=\bold{H_k}\ \Lambda_k +a_{k+1}r_{k+1}(1-\Lambda_{k}),
\end{equation*}
\noindent and \ $${\bf {H_{k+1}'}}={\bf {H_k'}}\Lambda_k
+a_{k+1}s_{k+1}(1-\Lambda_k)\bold{g_{k+1}}+({\bf
H_k}-a_{k+1}r_{k+1})\Lambda'_k.$$ Notice that, on the open set
$U_{k+1}$\,, we have that \ $\phi_{i,j}(g_i)\equiv 1 $, \ whenever
$i+j\le k$. Indeed, on the one hand, \  if $x\in U_{k+1}$, and
$i\in \{1,...,k\}$, then $g_i(x)<\gamma_{i,k+1-i}\le
\frac{\gamma_i+\gamma_{i,j}}{2}$, whenever $i+j\le k$.  \, On the
other hand, \ $\phi_{i,j}(t)\equiv 1$ \ if \ $t\le
\frac{\gamma_i+\gamma_{i,j}}{2}$. \ Therefore
$h_i|_{U_{k+1}}=\varphi_i(g_i)$, \ for every \ $i=1,....,k$, \ and
$$\Lambda_{k}=\frac{\sum_{i=1}^{k} \varphi_i(g_i)}{\sum_{i=1}^{k}
\varphi_i(g_i)+h_{k+1}}.$$ By computing $\Lambda'_k$ \ in
$U_{k+1}$, we obtain
$\Lambda'_k=\xi_{k,1}\bold{g_1}+...+\xi_{k,k+1} \bold{g_{k+1}}$,
where the coefficients $\xi_{k,1},...,\xi_{k,k+1}$ \ are
continuous functions of the following form:
 \begin{align*}\label{derivada de Lambda k}
&\xi_{k,j}= \frac{-\varphi_{k+1}(g_{k+1})\,\phi_{j,k+1-j}'(g_j)\,
(\prod_{i=1;\,i\not=j}^{k}\phi_{i,k+1-i}(g_i))(\sum_{i=1}^k
h_i)+h_{k+1}\varphi_j'(g_j)}{(\sum_{i=1}^{k+1}h_i)^2}, \quad
j=1,...,k \\
&\xi_{k,k+1}= \frac{-\varphi_{k+1}'(g_{k+1})\, ( \prod_{i=1}^k
\phi_{i,k+1-i}(g_i))\, (\sum_{i=1}^k h_i)}{(\sum_{i=1}^{k+1}
h_i)^2}.\,
\end{align*}

 Thus, if $x\in U_{k+1}\cap(U_1\cup ...\cup U_k)$, the
coefficients $\sigma_{k+1,1},....,\sigma_{k+1,k+1}$ \ for ${\bf
H_{k+1}'}$ have the following form,
\begin{align*}
&\sigma_{k+1,j}=\sigma_{k,j}\Lambda_k+({\bf H_k}-a_{k+1}r_{k+1})\xi_{k,j}, \qquad \text{ for } \ j=1,...,k\\
&\sigma_{k+1,k+1}=a_{k+1}s_{k+1}(1-\Lambda_k)+({\bf H_k}-a_{k+1}r_{k+1})\xi_{k,k+1}.\\
\end{align*}
\noindent Notice that  in $U_{k+1}\cap (U_1\cup ...\cup U_k)$, \
$a_{k+1}s_{k+1}>0$, \ $\Lambda_k>0$, \ $1-\Lambda_k>0$, \
$\xi_{k,j}\ge 0$, \ for every \ $j=1,...,k$, \
$\sum_{j=1}^k\xi_{k,j}> 0$ \ and \ $\xi_{k,k+1}<0$. Therefore, if
$H_k-a_{k+1}r_{k+1}\le 0$, the coefficient $\sigma_{k+1,k+1}>0$.
When $H_k-a_{k+1}r_{k+1} \ge 0$ and $\sigma_{k,j}>0$,
 the coefficient $\sigma_{k+1,j}>0$ (recall that, from the step $k$ we know that,
  for every $x\in U_1\cup ...\cup U_k$ there exists
 at least one $j\in \{1,...,k\}$ \ with \ $\sigma_{k,j}>0$). Hence, if
\ ${\bf H_{k+1}'}(x)= 0$ \ for some \ $x\in U_{k+1}\cap (U_1 \cup
...\cup U_k)$, \ there necessarily exists $\varrho\not=0$ such
that $D(x)=\varrho(\sigma_{k+1,1}(x)g_1
+...+\sigma_{k+1,k+1}(x)g_{k+1})$, \ that is
$D(x)\in[g_1,\,...,g_{k+1}].$ \

\medskip

\noindent In fact we can be more accurate and obtain that if \
${\bf H_{k+1}'}(x)=0$, \ $x\in U_{k+1}\cap (U_1\cup...\cup U_k)$
\and  $x\not\in \cup_{j\in F}U_j$ \ for some proper subset \
$F\subset\{1,...,k\}$, then \ $D(x)\in
\operatorname{span}\,\{g_j:\ j\in\{1,...,k+1\}\setminus F\}.$
Indeed, from step $k$ we know that,  if $x \in(U_1\cup...\cup U_k)
\setminus U_j$, where $j\in\{1,...,k\}$, then $\sigma_{k,j}(x)=0$.
Now, if $j\in F$ and $j=1$, it is clear that the functions
$\varphi_1'(g_1)$ and $\phi_{1,k}'(g_1)$ vanish outside $U_1$.
This implies $\xi_{k,1}(x)=0$ and consequently
$\sigma_{k+1,1}(x)=0$. If \ $j\in F$ and  $2\le j\le k$, since
$x\in U_{k+1}$ we know that
\begin{equation*}
g_1(x)<\gamma_{1,k}<\gamma_{1,j-1}\,,.......,g_{j-1}(x)<\gamma_{j-1,k+2-j}<\gamma_{j-1,1}\,,
\end{equation*}
and then necessarily \ $g_j(x)\le \gamma_j$. Since the functions \
$\varphi_j'(g_j)$ \ and \ $\phi_{j,k+1-j}'(g_j)$ \ vanish whenever
$g_j \le \gamma_j$, \ we deduce  \ $\xi_{k,j}(x)=0$ \ and thus \
$\sigma_{k+1,j}(x)=0$.

\medskip

\noindent Let us now define the sets
\begin{align*}
Z_{k+1,1}&=\begin{cases}\{z_{k+1}\}, & \text{ if } \ z_{k+1}\in
U_{k+1}\setminus (\overline{U}_1\cup ...\cup \overline{U}_k)\\
\emptyset, & \text{ otherwise }
\end{cases}\\
Z_{k+1,2}&=Z_{k+1}\cap U_{k+1}\cap (U_1\cup ...\cup U_k).
\end{align*}

Now, let us check that $Z_{k+1,2}\subset\cup_{g\in N_{k+1}}O_g$.
Indeed, if $x\in Z_{k+1,2}$\,, \ there are constants \ $1\le s\le
k$ \ and \ $1\le i_1<...<i_k\le k $\,, such that $x\in U_{k+1}\cap
U_{i_1}\cap ...\cap U_{i_s}$ \ and \ $x\not\in \cup_{j\in F} U_j$,
where $F=\{1,...,k\}\setminus \{i_1,...,i_s\}$. From the preceding
assertion, $D(x)\in [g_{i_1},...,g_{i_s},g_{k+1}]$. From the
definition of $M_{k+1,i_1,...,i_s}$ \ and the fact that
$U_{k+1}\subset U_{k+1}'$\,, \ we obtain that \ $x\in
M_{k+1,i_1,...i_s}\subset M_{k+1}$. \ Since $x\in U_{k+1}$,\ we
have that $g_1(x)<\gamma_{1,k}$\,,...,$g_k(x)<\gamma_{k,1}$. We
apply Fact \ref{k+1}(2) to conclude that there is $g\in N_{k+1}$
such that $x\in O_g$.

\medskip

In the case when $Z_{k+1,1}=\{z_{k+1}\}\not\in \cup_{g\in
N_{k+1}}\overline{O}_g$\,, we select, if necessary, a larger
$t_{k+1}$, \ with $t_{k+1}<l_{k+1}$\,,
 so that \ $z_{k+1}\not\in \cup_{g\in
N_{k+1}}\overline{B}_g$ \ for every \ $g\in N_{k+1}$. Since the
norm is LUR  and $D(z_{k+1})=g_{k+1}$ \ we may select numbers \
$0<t_{k+1}'<l_{k+1}'<1$ \ and open slices, which are neighborhoods
of \ $z_{k+1}$ \ defined by
\begin{equation*}
O_{g_{k+1}}:=\{x\in S: g_{k+1}(x)>l_{k+1}'\} \quad \text{ and }
\quad B_{g_{k+1}}:=\{x\in S: g_{k+1}(x)>t_{k+1}'\},
\end{equation*}
satisfying  \ $O_{g_{k+1}}\subset B_{g_{k+1}}\subset \{x\in S: \
g_1(x)<\gamma_{1,k}'\,,...,\ g_k(x)<\gamma_{k,1}'\,, \
g_{k+1}(x)>\gamma_{k+1}'\} $ \ and \ $\dist(B_{g_{k+1}},\,O_g)>0$,
\  for every $g\in N_{k+1}$.
 In this case, we define \ $\Gamma_{k+1}=N_{k+1}\cup\{g_{k+1}\}$.

\medskip

 Now, if $Z_{{k+1},1}=\{z_{k+1}\}\in \cup_{g\in N_{k+1}}{\overline{O}_g}$,\
 we select, if necessary a smaller constant $l_{k+1}$\,, with $0<t_{k+1}<l_{k+1}<1$\,, so that
 $Z_{k+1,1}=\{z_{k+1}\}\in \cup_{g\in N_{k+1}}{O_g} $. In this
 case, and also when
 $Z_{k+1,1}=\emptyset$, we  define $\Gamma_{k+1}=N_{k+1}$.

Notice that, in any of the cases mentioned above, Fact \ref{k+1}
clearly
 holds for the (possibly) newly selected real numbers $t_{k+1}$ and $l_{k+1}$.

\medskip

 Then, the distance between  any two sets  $B_{g}$, \  $B_{g'}$, \ where \ $g,g'\in
 \Gamma_1\cup... \cup \Gamma_{k+1}$, \ and \
 $g\not=g'$, \  is strictly positive. \
Moreover \ $Z_{{k+1},1}\cup Z_{{k+1},2}\subset \cup_{g\in
\Gamma_{k+1}} O_g \subset  \cup_{g\in \Gamma_{k+1}} B_g\subset
U_{k+1}'\subset R_{k+1}$. Therefore, $Z_{k+1}=Z_1\cup ...\cup Z_k
\cup Z_{{k+1},1}\cup Z_{{k+1},2}\subset \cup_{g\in \Gamma_1 \cup
... \cup \Gamma_{k+1} }O_g \subset \cup_{g\in \Gamma_1 \cup ...
\cup \Gamma_{k+1} }B_g \subset U_1\cup...\cup U_{k+1}=R_1\cup ...
\cup R_{k+1}$. Also, recall that \ $\dist(B_g, R_{k+1}^c)>0$, \
for every \ $g\in \Gamma_{k+1}$ \ and \ $\dist(B_g, (U_1\cup
...\cup U_{k+1})^c)>0$,\  for every \ $g\in \Gamma_1\cup ... \cup
\Gamma_{k+1}$.

\medskip

Finally, let us notice that, by combining the results obtained in
step k+1, we deduce that  \ ${\bf
H_{k+1}'}=a_{k+1}s_{k+1}\bold{g_{k+1}}$ \ in \ $U_{k+1} \setminus
(U_1\cup ...\cup U_k)$ \ and \ ${\bf H_{k+1}'}={\bf H_k'}$ on
$(U_1\cup...\cup U_k)\setminus U_{k+1}$, \ and in general \  ${\bf
H_{k+1}'}=\sigma_{k+1,1}\,\bold{g_1}+\cdots+\sigma_{k+1,k+1}\,\bold{g_{k+1}}$
\ on \ $U_1\cup...\cup U_{k+1}$, \ where \ $\sigma_{k+1,i}$ \ are
continuous functions on \ $U_1\cup...\cup U_{k+1}$, \ for \
$i=1,....,k+1$,
 \ and for every
 $x\in U_1 \cup \dots \cup U_{k+1}$ \ there is at least one $i\in
 \{1,...,k+1\}$  \ such that $\sigma_{k+1,i}(x)>0$. Furthermore,  \
$\sigma_{k+1,j}(x)=0$  \ whenever \ $x\in (U_1 \cup ...\cup
U_{k+1}) \setminus U_j$,  \  \ $j\in \{1,...,k+1\}$.

\bigskip

Once we have defined, by induction, the functions $h_k$, \ $r_k$ \
and the constants $a_k$, for all $k\in \mathbb N$, we define
\begin{align*}
&H:S^+\longrightarrow \mathbb R\\  &H=\frac{\sum_{k=1}^\infty
a_kr_kh_k}{\sum_{k=1}^\infty h_k}.
\end{align*}
It is straightforward to verify that the family $\{U_k\}_{k\in
\mathbb N}$ of open sets of $S^+$ is a locally finite open
covering of $S^+$. Thus, for every $x\in S^+$ there is $k_x\in
\mathbb N$ and  a (relatively  open in $S^+$) neighborhood  \
$V_x\subset S^+$ \ of \ $x$, such that \ $V_x\cap(\cup_{k>k_x}
U_k)=\emptyset$ \ and therefore \ $H|_{V_x}={\bf H_{k_x}}|_{V_x}$.
Thus $H$ is $C^p$ smooth whenever the functions $\{h_k\}_{k\in
\mathbb N}$ are $C^p$ smooth.

\medskip

\begin{fact}
The function \ $H$ \ $3\,\varepsilon$-approximates \ $F$ \ in
$S^+$.
\end{fact}
\noindent {\em Proof.} Recall that the oscillation of $F$ in $U_k$
is less that $\varepsilon$, \ the oscillation of $r_k$  in $U_k$
is less than $\frac{\varepsilon}{\,|a_k|\,}$\, ,  \
$|a_k-F(x_k)|<\varepsilon$ \ and \ $r_k(x_k)=1$,\  for every $k\in
\mathbb N$. Now, if $h_k(x)\not=0$, then $x\in U_k$ and
\begin{align}\label{akrk-F}
|a_kr_k(x)-F(x)|& \le|a_kr_k(x)-a_kr_k(x_k)|+|a_kr_k(x_k)-F(x)|\\
\notag &= |a_k||r_k(x)-r_k(x_k)|+|a_k-F(x)|\\ \notag &\le
|a_k|\frac{\varepsilon}{|a_k|}+|a_k-F(x_k)|+|F(x_k)-F(x)| \le
3\varepsilon.
\end{align}
Hence,
\begin{align*}
|H(x)-F(x)|=\frac{\bigl|\sum_{k=1}^\infty(a_kr_k(x)-F(x))\,h_k(x)
\bigr|}{\sum_{k=1}^\infty h_k (x)}\le
\frac{\sum_{k=1}^\infty|a_kr_k(x)-F(x)|\,h_k(x)
}{\sum_{k=1}^\infty h_k (x)}\le 3\varepsilon. \Box
\end{align*}

Let us denote by \ $C$ \ the critical points of $H$ in $S^+$.
Since for every $x\in S^+$, there is $k_x\in \mathbb N$ and a
(relatively open in $S^+$) neighborhood $V_x\subset S^+$ of $x$
such that $V_x\cap(\cup_{k>k_x} U_k)=\emptyset$, we have that
$H|_{V_x}=\bold{H_{k_x}}|_{V_x}$ \ and $C\subset \cup_k Z_k$.
Recall that $\cup_k Z_k \subset \bigcup \{O_g: \ g\in
\cup_k\Gamma_{k}\} \subset \bigcup \{B_g:\ g\in \cup_k\Gamma_k\}$,
\ the oscillation of $F$ on $B_g$ is less than $\varepsilon$ \ and
$\dist(B_g,B_{g'})>0$, \ for every $g,g'\in \cup_k \Gamma_k$ with
$g\not=g'$. Furthermore, from the inductive construction of the
sets $\{B_g: \ g\in\cup_k \Gamma_k \}$, \ it is straightforward to
verify that (i) for every $k>1$, if \ $g\in \Gamma_k $ \ and \
$g'\in \cup_{m>k}\Gamma_m$, \ then $\dist(B_g,\,B_{g'})\ge
\gamma_k'-\gamma_{k,1}'>0$ \ and  (ii) if $g\in \Gamma_1$ \ and \
 $g'\in \cup_{m>1}\Gamma_m$, then $\dist(B_g,\,B_{g'})\ge
t_1-\gamma_{1,1}'>0$. Therefore, for every $g\in\cup_k\Gamma_k$,
\begin{equation}\label{distanciasBs} \dist(B_g\,, \,\bigcup\{B_{g'}:\ g'\in \cup_k\Gamma_k,\ \
g'\not=g\})>0.
\end{equation} We relabel the countable families of open slices \
$\{O_g\}_{g\in \cup_k\Gamma_k}$ \ and \ $\{B_g\}_{g\in
\cup_k\Gamma_k}$ \ as \ $\{O_n\}$, \ $ \{B_n\}$, \ respectively.
Notice that the set $\cup_n\overline{B}_n$ is a (relatively)
closed set in $S^+$. Indeed, if $\{x_j\}\subset
\cup_n\overline{B}_n$ and $\lim_j x_j=x\in  S^+$, since $\cup_n
U_n'$ \ is also  a locally finite open covering of $S^+$, there is
$n_x$ and a (relatively open in $S^+$) neighborhood $W_x\subset
S^+$ of $x$\,, \ such that $W_x\cap(\cup_{n> n_x}U_n')=\emptyset$.
In addition, from the construction of the family $\{B_n\}$, \
there is \ $N\in\mathbb N$ \ such that \ $\cup_{n>N}
\overline{B}_n\subset \cup_{n> n_x}U_n'$\,, and thus there is \
$j_0\in \mathbb N$ \ with  \ $\{x_{j}\}_{j> j_0}\subset
\cup_{n=1}^{N} \overline{B}_n$. Hence \ $x\in \cup_{n=1}^{N}
\overline{B}_n\subset\cup_{n} \overline{B}_n$.

\bigskip

 Let us denote \ $\mathcal{B}_n=\Phi^{-1}(B_n)$ \ and \
$\mathcal{O}_n=\Phi^{-1}(O_n)$, \ for every $n\in \mathbb N$.
\begin{fact}
$\mathcal{O}_n$ and $\mathcal{B}_n$ are open, convex and bounded
subsets of $X$, for every $n\in\mathbb N$.
\end{fact}
\noindent {\em Proof.}  Since $\Phi$ is continuous, it is clear
that $\mathcal{O}_n$ \ and \ $\mathcal{B}_n$ are open sets.
 The sets ${O}_n$ and ${B}_n$ are slices of the
 form $R=\{x\in S: \ b(x)>\delta\}$ \ for some \ $b\in S^*$ \ and
 $\delta>0$ such that $\dist(R,X\times\{0\})>0$. Let us prove that
 $\mathcal{R}:=\Phi^{-1}(R)$ is  convex and bounded in $X$. First, let us check that the cone
 in $Y$ generated by $R$ \ and  defined by $$\co(R)=\{\lambda x:\ x\in R, \
 \lambda>
 0\}=\{x\in Y: b(\textstyle{\frac{x}{\,|x|\,}})>\delta\}$$ is a convex set: consider $0\le \alpha \le 1$ and
 $x,x'\in \co(R)$. Then, \begin{align*} b(\alpha x+(1-\alpha)x')&=\alpha b(x)+
 (1-\alpha )b(x')>
 \alpha\delta|x|+(1-\alpha)\delta|x'|\\&=\delta|\alpha x|+\delta
 |(1-\alpha) x'|\ge \delta |\alpha x+
 (1-\alpha) x'|,\end{align*}
and this implies that $\alpha \,x+(1-\alpha)\,x'\in \co(R)$.
Therefore, the intersection of the two convex sets \ $\co(R)\cap
(X\times\{1\})=\Pi^{-1}(R)$ \ is convex. Now, it is clear that
$\mathcal{R}=\Phi^{-1}(R)=i^{-1}(\Pi^{-1}(R))$ is convex as well.

Let us prove that \ $\Pi^{-1}(R)$ \ is bounded in \ $Y$. Consider
the linear bounded operator \ $\pi_2:Y=X\oplus \mathbb
R\longrightarrow \mathbb R$,\ $\pi_2(x,r)=r$\,, \ for every
$(x,r)\in X\oplus\mathbb R$. Then, \
$\Pi^{-1}(y)=\frac{y}{\pi_2(y)}$\,, for every $y\in S^+$. On the
one hand,  \ $d:=\dist(R,\,X\times\{0\})>0$ \ and then \
$$\ \ \ \pi_2(x,r)=r=\frac{|(x,r)-(x,0)|}{|(0,1)|}\ge
\frac{d}{|(0,1)|}:=s>0, \quad \text{ for every } \ (x,r)\in R. $$
On the other hand,
\begin{equation*}
\bigr|\Pi^{-1}(y)\bigl|=\frac{|y|}{\pi_2(y)}=\frac{1}{\pi_2(y)}\le
\frac{1}{s}, \quad \text{ for every  } y\in R,
\end{equation*}
and thus $\Pi^{-1}(R)$ is bounded. Since the norm \ $||\cdot||$ \
considered on \ $X\times\{0\}$ \ (defined as \ $||(x,0)||=||x||$)
\ and the restriction of the norm \ $|\cdot|$ \ to \
$X\times\{0\}$ are equivalent norms on \ $X\times\{0\}$, there
exist constants $m,M>0$ such that
$m\,||x-x'||\le|i(x)-i(x')|=|(x,1)-(x',1)|=|(x-x',0)|\le M
||x-x'||$, \ for every $x,x'\in X$. Hence,
\begin{align*}\hspace{1cm}||\Phi^{-1}(y)||&=||i^{-1}(\Pi^{-1}(y))-i^{-1}(0,1)||\le
\frac{1}{m}|\Pi^{-1}(y)-(0,1)| \\& \le
\frac{1}{m}\,\bigl(|\Pi^{-1}(y)|+|(0,1)|\bigr) \le \frac{1
+s|(0,1)|}{s\cdot m},\end{align*} for every $y\in R$, what shows
that \ $\mathcal{R}$ \ is bounded in \ $X$. $\Box$

\medskip

\begin{fact}
$\overline{\mathcal{O}}_n $ \ and \ $\overline{\mathcal{B}}_n$ \
are (closed convex and bounded) $C^p$ smooth bodies, \ for every
$n\in \mathbb N$.
\end{fact}
\noindent {\em Proof.}  We already know that these sets are
closed, convex and bounded bodies, hence it is enough to prove
that their boundaries $\partial\mathcal{O}_n$ and $\partial
\mathcal{B}_n$ are $C^p$ smooth one-codimensional submanifolds of
$X$. Since $\partial\mathcal{B}_n=\Phi^{-1}(\partial B_n)$,
$\partial\mathcal{O}_n=\Phi^{-1}(\partial O_n)$, and $\Phi$ is a
$C^p$ diffeomorphism, this is the same as showing that $\partial
O_n$ and $\partial B_n$ are $C^p$ smooth one-codimensional
submanifolds of $S$. But, if $O_n$ is defined by $O_n=\{y\in S:
g_{n}(y)>\beta_{n}\}$, we have that $\partial O_n$ is the
intersection of $S$ with the hyperplane $X_n=\{y\in Y:
g_{n}(y)=\beta_{n}\}$ of $Y$, and $X_n$ is transversal to $S$ at
every point of $\partial O_n$ (otherwise the hyperplane $X_n$
would be tangent to $S$ at some point of $\partial O_n$ and, by
strict convexity of $S$, this implies that $\partial O_n=X_n\cap
S$ is a singleton, which contradicts the fact that $O_n$ is a
nonempty open slice of $S$), hence the intersection $\partial
O_n=S\cap X_n$ is a one-codimensional submanifold of $S$. The same
argument applies to $\partial B_n$. $\Box$

\medskip

\begin{fact} $\textrm{dist}(\mathcal{O}_n, X\setminus
\mathcal{B}_n)>0$  \ and
$\dist\left(\mathcal{B}_n,\cup_{m\not=n}\mathcal{B}_{m}\right)>0$,
\ for every $n\in\mathbb N$.
\end{fact}
\noindent {\em Proof.} This is a consequence of the fact that \
$\dist(O_n, S^+\setminus B_n)>0$,   \
$\dist\left(B_n,\cup_{m\not=n}B_{m}\right)>0$, \ and \ $\Phi$ is
Lipschitz. Indeed,  on the one hand, recall that
$|i(x)-i(x')|=|(x-x',0)|\le M||x-x'||$, \ for every $x,x'\in X$.
On the other hand,
\begin{equation*}\hspace{1cm}\bigl|\Pi(y)-\Pi(y')\bigl|=\frac{|\,y\,|y'|-y'\,|y|\,|}{|y|\,|y'|}=
\frac{\bigr|y\,(|y'|-|y|)+(y-y')|y|\bigl|}{|y|\,|y'|}\le
\frac{2|y-y'|}{|y'|}\le \frac{2}{\zeta}|y-y'|\,,\end{equation*}
for every $y,y'\in X\times \{1\}$, \ where
$\zeta=\dist(0,X\times\{1\})>0.$ Therefore, $|\Phi(x)-\Phi(x')|\le
\frac{2M}{\zeta}||x-x'||$, \ for every $x,x'\in X$. Now, if two
sets $A,A'\subset S^+$ satisfy that $\dist(A,A')>0$, then
$\dist(A,A')\le |a-a'|\le \frac{2M}{\zeta}||\Phi^{-1}(a)-\Phi
^{-1}(a')||$, for every $a\in A$,\  $a'\in A'$. Therefore,
$0<\dist(A,A')\le \frac{2M}{\zeta}\dist(\Phi^{-1}(A),\,\Phi
^{-1}(A'))$. $\Box$

\medskip

\begin{fact}
For every $n\in\mathbb N$, there exists a $C^p$ diffeomorphism \
$\Psi_n$ \ from \ $X$ \ onto \ $X\setminus
\overline{\mathcal{O}}_n$ \ such that \ $\Psi_n$  \ is the
identity off \ $\mathcal{B}_n$.
\end{fact}
\noindent {\em Proof.} Assume that $0\in \mathcal{O}_n$. Since
$\dist(\mathcal{O}_n, X\setminus \mathcal{B}_n)>0$, \ there is \
$\delta_n>0$ \ such that \
$\dist((1+\delta_n)\mathcal{O}_n\,,\,\mathcal{B}_n)>0$.\ We can
easily construct a $C^{p}$ smooth radial diffeomorphism \
$\Psi_{n,2}$ \ from \ $X\setminus\{0\}$ \ onto \ $X\setminus
\overline{\mathcal{O}}_n$ \ satisfying \ $\Psi_{n,2}(x)=x$ \ if \
$x\notin (1+\delta_n)\mathcal{O}_n$. Indeed, take a \ $C^{\infty}$
\ smooth function $\lambda_n:[0,\infty)\longrightarrow [1,\infty)$
satisfying that $\lambda_n(t)=t$ for $t\geq 1+\delta_n$,
$\lambda_n(0)=1$ \ and \ $\lambda'_n(t)>0$ for $t>0$, and define
$$
\Psi_{n,2}(x)=\lambda_n(\mu_n(x))\,\frac{x}{\mu_n(x)},
$$
for $x\in X\setminus\{0\}$, where $\mu_n$ is the Minkowski
functional of $\overline{\mathcal{O}}_n$, which is $C^p$  smooth
on $X\setminus\{0\}$.

Now, since $0\in \mathcal{O}_n$, there is $\alpha_n>0$ such that
$\alpha_n B_{||\cdot||}\subset \mathcal{O}_n$\,. According to
\cite[Proposition 3.1]{Dobrowolski1} and \cite[Lemma
2]{Dobrowolski2} (see also \cite{Az}), there exists a \ $C^p$ \
diffeomorphism \ $\Psi_{n,1}$ \ from $X$ \ onto \
$X\setminus\{0\}$ \  such that $\Psi_{n,1}$ is the identity off
$\alpha_n B_{||\cdot||}$ (this set may be regarded as the unit
ball of a equivalent $C^p$ smooth norm on $X$).

Then, the composition \ $\Psi_{n}:=\Psi_{n,2}\circ\Psi_{n,1}$
 \ is a \ $C^p$ diffeomorphism
from \ $X$ \ onto  \ $X\setminus\overline{\mathcal{O}}_n$ \ such
that \ $\Psi_{n}$ \ is the identity off \ $\mathcal{B}_n$. If
$0\not\in \mathcal{O}_n$, select $\omega_n\in \mathcal{O}_n$ \ and
repeat the above construction of the diffeomorphism with the sets
\ $\mathcal{O}_n-\omega_n$ \ and $\mathcal{B}_n-\omega_n$. Then,
$\Psi_{n}:=\tau_{\omega_n}\circ\Psi_{n,2}\circ\Psi_{n,1}\circ\tau_{-\omega_n}$
is the required $C^p $ diffeomorphism,\ where \
$\tau_\omega(x)=x+\omega$. $\Box$

\medskip

Now, the infinite composition $\Psi=\bigcirc_{n=1}^\infty \Psi_n$
is a well-defined $C^p$ diffeomorphism from $X$ onto $X\setminus
\cup_n \overline{ \mathcal{O}}_n$\,, which is the identity outside
$\cup_n \mathcal{B}_n$ and $\Psi(\mathcal{B}_n)\subset
\mathcal{B}_n$. This follows from the fact that, for every $x\in
X$, there is an open neighborhood $V_x$ and $n_x\in\mathbb N$ such
that $V_x\cap(\cup_{n\not=n_x} \mathcal{B}_n)=\emptyset$, and
therefore $\Psi|_{V_x}= \Psi_{n_x}|_{V_x}$.

\medskip

 Finally, let us check that the $C^p$ smooth function \begin{align*}
 &g:X\longrightarrow \mathbb
 R \\ &g:=H\circ \Phi\circ\Psi\end{align*}$4\varepsilon$-approximates $f$
 on $X$  and $g$ does not have critical points. Indeed, for every
 $x\in X$, if $\Psi(x)\not=x$, then there is $\mathcal{B}_{n_x}$ such that $x\in
 \mathcal{B}_{n_x}$. Since the oscillation of $f$ in $\mathcal{B}_{n_x}$ is less than $\varepsilon$
 and $\Psi(x)\in \mathcal{B}_{n_x}$, we can deduce that
 $|f(\Psi(x))-f(x)|<\varepsilon$, for every $x\in X$. Recall that $F\circ \Phi=f$ \ and \
 $|H(x)-F(x)|<3\varepsilon$, for every $x\in S^+$. Then,
\begin{align}\label{g-f}
|g(x)-f(x)|&=|H\circ\Phi(\Psi(x))-F\circ\Phi(x)|\\ \notag &\le
|H(\Phi(\Psi(x)))-F(\Phi(\Psi(x)))|+
|F\circ\Phi(\Psi(x))-F\circ\Phi(x)|\\ \notag &\le 3\varepsilon
+\varepsilon=4\varepsilon,
\end{align}
for every $x\in X$. Since $\Phi$ and $\Psi$ are $C^p$
diffeomorphisms, we have that $g'(x)=0$ if and only if
$H'(\Phi(\Psi(x)))=0$. \ For every $x\in X$, \ $\Psi(x)\not\in
\cup_n\overline{\mathcal{O}}_n$ \ and thus \ $\Phi(\Psi(x))\not\in
\cup_n \overline{O}_n$. It follows that \
$H'(\Phi(\Psi(x)))\not=0$ \ and \ $g$ does not have any critical
point.

\medskip

Before finishing the proof, let us say what additional precautions
are required in the case when $\varepsilon$ is a strictly positive
continuous function:
\begin{itemize} \item the slices \ $S_k=\{x\in S: \
f_k(x)>\delta_k\}$, \ $(k\in\mathbb N)$ \ are selected with the
additional property that the oscillation of the two functions \
$F$ \ and \ $\overline{\varepsilon}=\varepsilon \circ \Phi^{-1}$ \
in \ $S_k$ \ are less than \
$\frac{\overline{\varepsilon}(y_k)}{2}$, \ where $y_k$ is the
point of $S^+ $ satisfying $f_k(y_k)=1$; \ this implies, in
particular,  that \
$\frac{1}{2}\,\overline{\varepsilon}(y_k)<\overline{\varepsilon}(x)<\frac{3}{2}\,\overline{\varepsilon}(y_k)$,\
for every $x\in S_k$; \item the real numbers \ $a_k\in \mathbb
R^*$ \ satisfy that \
$|a_k-F(x_k)|<\frac{\overline{\varepsilon}(y_k)}{2}$; \item the
oscillation of $r_k$ in $S_k$ is less than
$\frac{\overline{\varepsilon}(y_k)}{|a_k|}$.
\end{itemize} From the above conditions and  inequality  \eqref{akrk-F}, it can be deduced that
if \  $x\in U_k$\,, then $|a_kr_k(x)-F(x)|\le
2\,\overline{\varepsilon}(y_k)<4\overline{\varepsilon}(x).$ \ From
this, it can be obtained that $|H(x)-F(x)|\le
4\overline{\varepsilon}(x)$, \ for every \ $x\in S^+$.
Equivalently, $|H\circ \Phi(x)-F\circ \Phi(x)|=|H\circ
\Phi(x)-f(x)|<4\,\varepsilon(x)$, \ for every $x\in X$. Now, if
$x\not=\Psi(x)$, then there is $\mathcal{B}_{n_x}$ such that
$x,\Psi(x)\in \mathcal{B}_{n_x}$. Thus, $|f(\Psi(x))-f(x)|<
\frac{\varepsilon(\Phi^{-1}(y_{n_x}))}{2}<\varepsilon(x)$. Now,
from inequality \eqref{g-f}, we obtain: (a) if $x\in B_{n_x}$ \
for some $n_x$,\ then \  $|g(x)-f(x)|\le
4\,\varepsilon(\Psi(x))+\varepsilon(x)\le
6\,\varepsilon(\Phi^{-1}(y_{n_x}))+\varepsilon(x)\le
13\,\varepsilon(x)$, and \ (b) if \ $x\not\in\cup_nB_n$\,, then
$|g(x)-f(x)|\le 4\,\varepsilon(\Psi(x))=4\,\varepsilon(x).$ \ This
finishes the proof of Theorem \ref{approximation theorem}. $\Box$

\bigskip

\begin{rem}  {\rm  The construction of the function $g$ with no critical points
that approximates $f$ with a constant $\varepsilon>0$, is
considerably shorter in the case that either (i) $X=\ell_2(\mathbb
N)$  (and we use West Theorem \cite{West})\ or \ (ii) $X$ is
non-reflexive and the norm $|\cdot|$ considered on $Y$ can be
constructed with the additional property that the set $\{ f\in
Y^{*}: f \textrm{ does not attain its norm} \}$ contains a dense
subspace (except the zero functional).

Indeed, in the first case, we can define as $|\cdot|$ the standard
norm on $\ell_2(\mathbb N)$.
 In both cases, the use of the auxiliary
functions $r_n$\, is not required,\  we can consider the slice \
$R_n:=S_n$ (that is, the additional construction of the sequence
of slices $\{R_n\}$ is not required) \ and \ we can select for
every $n\in \mathbb N$, \ any strictly decreasing sequence \
$\{\gamma_{n,i}\}_i$ \ such that \
$\lim_i\gamma_{n,i}=\delta_n$\,. \ Then, \ let us choose a
non-zero functional $w\in Y^*\setminus [f_n: n\in \mathbb N]$
(where $[f_n:n\in \mathbb N]$ denotes the space of all finite
linear combinations of the set $\{f_n:\, n\in \mathbb N\})$ with
$|w|^*<\varepsilon$, \ and define \
$H=\frac{\sum_ia_ih_i}{\sum_ih_i}+w$ \ and \
$\bold{H_n}=\frac{\sum_{i=1}^na_ih_i}{\sum_{i=1}^nh_i}+w$, \ for
every $n\in \mathbb N$. We obtain  in the case (i),  that \ $Z_n$
\ (the critical points of $\bold{H_n}$), is included in the
compact set $D^{-1}([f_1,...,f_n,w]\cap S^*)$. Therefore, the set
$C$ of critical points of $H$ and thus the set $\mathcal{C}$  of
critical points of the composition \ $ H \circ \Phi $, \ are
closed and locally compact sets of $S^+$ and $\ell_2$\,,
respectively.  \ Now, $g$ is obtained, by applying West Theorem
\cite{West}, considering a $C^\infty$ \ deleting diffeomorphism
$\Psi$ from $\ell_2(\mathbb N)$ onto $\ell_2(\mathbb N)\setminus
\mathcal{C}$, with the additional property that the family \
$\{(x,\,\Psi(x)):\ x\in\ell_2(\mathbb N)\}$ \ refines the open
convering  $ \{\Phi^{-1}(S_{n});\ n\in \mathbb N\}) $. Finally, \
we can define \ $g:= H \circ \Phi \circ \Psi$.

\noindent In the case (ii), we can select the family
$\mathcal{G}=\{f_n:\ n\in \mathbb N\}\cup\{w\}$ \ with the
additional requirement that \ $[\mathcal{G}]\setminus\{0\}$ \ is
included in the set of non-norm attaining functionals. \ Thus, the
sets of critical points of both $\bold{H_n}$ and $H$ are empty.
Therefore, the set of critical points of $g:=H\circ\Phi$ is empty
and $g$ approximates $f$. Notice that this case is particularly
interesting because the use of a deleting diffeomorphism is not
required.
}
\end{rem}

\begin{center}
{\bf Acknowledgements}
\end{center}
\noindent  This research was carried out during
Jim\'{e}nez-Sevilla's stay at the Mathematics Department of Ohio
State University; Jim\'{e}nez-Sevilla wishes to thank very
specially Peter March and Boris Mityagin for their kind
hospitality. Azagra thanks Gilles Godefroy and Yves Raynaud for
all their help during his stay at Institut de Math\'{e}matiques de
Jussieu (Universit\'{e} Paris 6).

\vspace{3mm} \noindent Departamento de An\'{a}lisis
Matem\'{a}tico. Facultad de Ciencias Matem\'{a}ticas. Universidad
Complutense. 28040 Madrid, SPAIN. \\ \noindent {\em E-mail
addresses:} azagra@mat.ucm.es, mm\_jimenez@mat.ucm.es


\bigskip


\begin{thebibliography}{21}
\addcontentsline{toc}{chapter}{Bibliography}


\bibitem{Az}D. Azagra, {\em Diffeomorphisms between spheres and hyperplanes between
infinite-dimensional Banach spaces,  } Studia Math. {\bf 125}(2)
(1997), 179-186.

\bibitem{AC}
D. Azagra and M. Cepedello Boiso, {\em Uniform approximation of
continuous mappings by smooth mappings with no critical points on
Hilbert manifolds}, Duke Math. J. {\bf 124} (2004), 207--226.


\bibitem{AFGJL} D. Azagra, R. Fry, J. G\'{o}mez Gil, J.A. Jaramillo,
M. Lovo, {\em $C^1$-fine appoximation of functions on Banach
spaces with unconditional Basis}, Quart. J. Math. {\bf 56} 2005,
13-20.

\bibitem{AJ}
D. Azagra and M. Jim\'{e}nez-Sevilla, {\em The failure of Rolle's
Theorem in infinite dimensional Banach spaces}, J. Funct. Anal.
{\bf 182} (2001), 207--226.


\bibitem{Bates1}
S. M. Bates, {\em On the image size of singular maps. I.}, Proc.
Amer. Math. Soc. {\bf 114} (1992), no. 3, 699-705.

\bibitem{Bates2}
S. M. Bates, {\em On the image size of singular maps. II.}, Duke
Math. J. {\bf 68} (1992), no. 3, 463-476.

\bibitem{Bates3}
S. M. Bates, {\em Toward a precise smoothness hypothesis in Sard's
theorem}, Proc. Amer. Math. Soc. {\bf 117} (1993), no. 1, 279-283.

\bibitem{Bates4}
S. M. Bates, {\em On smooth rank-1 mappings of Banach spaces onto
the plane}, J. Differential Geom. {\bf 37} (1993), no. 3, 729-733.

\bibitem{Bates-Moreira}
S. M. Bates and C. G. Moreira, {\em De nouvelles perspectives sur
le th\'{e}or\`{e}me de Morse-Sard}, C.R. Acad. Sci. Paris, t. {\bf
332}, S\'{e}rie I (2001), p. 13-17.

\bibitem{Be}
C. Bessaga, {\em Every infinite-dimensional Hilbert space is
diffeomorphic with its unit sphere},
 Bull. Acad. Polon. Sci., S\'{e}r. Sci. Math.
Astr. et Phys. {\bf 14} (1966), pp. 27-31.

\bibitem{DGZ}
R. Deville, G. Godefroy, and V. Zizler, {\em Smoothness and
renormings in Banach spaces}, vol. {\bf 64}, Pitman Monographs and
Surveys in Pure and Applied Mathematics,  Longman Scientific \&
Technical, Harlow, 1993.

\bibitem{Dobrowolski1}  T. Dobrowolski, {\em Smooth and R-Analitic
negligibility of subsets and extensions of homeomorphisms in
Banach spaces}, Studia Math. {\bf 65} (1979), 115--139.


\bibitem{Dobrowolski2}  T. Dobrowolski, {\em Relative classification
of smooth convex bodies}, Bull. Acad. Polon. Sci., S\'{e}r. Sci.
Math. {\bf 25} (1977), 309--312.


\bibitem{EellsMcAlpin}
J. Eells and J. McAlpin, {\em An approximate Morse-Sard theorem},
J. Math. Mech. {\bf 17} (1967/1968), 1055-1064.

\bibitem{Ferrer}
J. Ferrer, {\em Rolle's theorem fails in $\ell_{2}$}, Am. Math.
Monthly, vol. {\bf 103}, n. 2 (1996), 161-165.


\bibitem{Hajek}
P. H\'{a}jek, {\em Smooth functions on $c_{0}$}, Israel J. Math.,
{\bf 104} (1998), 17-27.

\bibitem{HJ} P. Hajek and M. Johanis, {\em Smooth approximations without
critical points}, Cent. Eur. J. Math. {\bf 1} (2003), no. 3,
284--291

\bibitem{Hirsch}
M. W. Hirsch, {\em Differential topology.} Graduate Texts in
Mathematics, No. 33. Springer-Verlag, New York-Heidelberg, 1976.

\bibitem{Kupka}
I. Kupka, {\em Counterexample to the Morse-Sard theorem in the
case of infinite-dimensional manifolds}, Proc. Amer. Math. Soc.
{\bf 16} (1965), 954-957.

\bibitem{Moreira}
C. G. Moreira, {\em Hausdorff measures and the Morse-Sard
theorem}, Publ. Mat. {\bf 45} (2001), 149-162.

\bibitem{Morse}
A. Morse, {\em The behavior of a function on its critical set},
Annals of Math. {\bf 40} (1939), 62-70.

\bibitem{Sard1}
A. Sard, {\em The measure of the critical values of differentiable
maps}, Bull. Amer. Math. Soc. {\bf 48} (1942), 883-890.

\bibitem{Sard2}
A. Sard, {\em Images of critical sets}, Annals of Math. {\bf 68}
(1958), 247-259.

\bibitem{Sard3}
A. Sard, {\em Hausdorff measure of critical images on Banach
manifolds}, Amer. J. Math. {\bf 87} (1965), 158-174.

\bibitem{Shk}
S. A. Shkarin, {\em On Rolle's theorem in infinite-dimensional
Banach spaces}, translation from Matematicheskie Zametki, vol.
{\bf 51}, no.3, pp. 128-136, March, 1992.

\bibitem{Smale}
S. Smale, {\em An infinite dimensional version of Sard's theorem},
Am. J. Math. {\bf 87} (1965), 861-866.

\bibitem{YomdinComte}
Y. Yomdin and G. Comte,  {\em  Tame geometry with applications in
smooth analysis}, Lecture Notes in Mathematics vol. {\bf 1834}
Springer-Verlag, 2004.

\bibitem{West} James E. West, {\em The diffeomorphic excision of closed local compacta
from infinite-dimensional Hilbert manifolds}, Compositio Math.
{\bf 21} (1969), 271-291.


\end{thebibliography}
\end{document}